\documentclass[preprint,12pt,authoryear]{elsarticle}
\usepackage[left=2.8cm,top=2.8cm,right=2.8cm,bottom=2.8cm]{geometry}
\usepackage{amsmath}
\usepackage{graphicx}
\usepackage{subcaption}
\usepackage[utf8]{inputenc}
\usepackage{times} 
\usepackage{latexsym,amsfonts,amssymb,amsmath,amscd,mathrsfs}
\usepackage{xargs}   
\usepackage{xspace}
\usepackage[dvipsnames]{xcolor}
\usepackage{float}
\usepackage{multirow}
\usepackage{graphicx,tikz,times,pgffor}
\usepackage{pgfplots}
\pgfplotsset{compat=1.17} 
\usetikzlibrary{calc,through,backgrounds,patterns}
\usepackage{tabulary}
\usepackage{tabularx}
\usepackage{changepage}
\usepackage{soul}
\usepackage{todonotes}
\usepackage[ruled,vlined]{algorithm2e}
\usepackage{hyperref}
\usepackage{setspace}
\usepackage{lineno}

\newcommand{\x} {\mathrm{x}}

\newcommand{\setX}{\mathcal{X}}
         
\newcommand{\X}{\mathbb{X}}

\newcommand{\setE}{\mathcal{E}}

\newcommand{\R}{\mathcal{R}}
       
\newcommand{\oR}{\overline{\mathcal{R}}}
    
\newcommand{\uR}{\underline{\mathcal{R}}}

\newtheorem {propt}{Property}[section]

\newtheorem {theo}{Theorem}[section]
\newtheorem {defn}{Definition}[section]
\newtheorem {cor}{Corollary}[section]
\newtheorem {lem}{Lemma}[section]
\newtheorem {rem} {Remark}[section]
\newtheorem {exmp} {Example}[section]

\newenvironment{pf}{\noindent \textbf{Proof:}}{}

\journal{}

\begin{document}
	
	\begin{frontmatter}
		
		\title{Minimizing the Epidemic Final Size while Containing the Infected Peak Prevalence in SIR Systems}

		\author[inst1]{J. Sereno}
		
		\affiliation[inst1]{organization={Institute of Technological Development for the Chemical Industry (INTEC), CONICET-Universidad Nacional del Litoral (UNL)},
			addressline={Guemes 3450}, 
			city={Santa Fe},
			postcode={3000}, 
			country={Argentina}}
		
		\author[inst1]{A. L. Anderson}
		
		\author[inst2]{A. Ferramosca}
		
		\affiliation[inst2]{organization={Department of Management, Information and Production Engineering, University of Bergamo},
			addressline={Via Marconi 5}, 
			city={Dalmine},
			postcode={24044}, 
			country={Italy}}
		
		\author[inst3,FIAS,shangai]{E.A. Hernandez-Vargas\corref{cor1}}
		\ead{esteban@im.unam.mx}
		
		\affiliation[inst3]{organization={Instituto de Matemáticas, UNAM},
			addressline={Boulevard Juriquilla 3001}, 
			city={Querétaro},
			postcode={76230}, 
			country={Mexico}}
		
		\affiliation[FIAS]{organization={Frankfurt Institute for Advanced Studies},  
			addressline={Ruth-Moufang-Str. 1, 60438},
			city={Frankfurt am Main}, country={Germany}}
		
		\affiliation[shangai]{organization={Director of 111 Oversea Talent Base for Intelligent Measurement \& Control for Complex Networked Systems and Applications}, city={Shanghai University}, country={China}}
		
		\author[inst1]{A. H. Gonz\'alez\corref{cor1}}
		\ead{alejgon@santafe-conicet.gov.ar}
		
		\cortext[cor1]{Corresponding author}
		
		\doublespacing

\begin{abstract}                          
Mathematical models are instrumental to forecast the spread of pathogens and to evaluate the effectiveness of non-pharmaceutical measures. A plethora of optimal strategies has been recently developed to minimize either the infected peak prevalence (IPP) or the epidemic final size (EFS). While most of the control strategies optimize a simple cost function along a fixed finite-time horizon, no consensus has been reached about how to simultaneously handle the IPP, the EFS, and the avoiding of new cycles of infections rebounding.
In this work, based on the characterization of the dynamical behaviour of SIR-type models under control actions (including the stability of equilibrium sets, in terms of the herd immunity), it is studied how to minimize the EFS while keeping - at any time - the IPP controlled. A procedure is proposed to tailor non-pharmaceutical interventions by separating transient from stationary control objectives and the potential benefits of the strategy are illustrated by a detailed analysis and simulation results related to the COVID-19 pandemic.
\end{abstract}


\end{frontmatter}


\section{Introduction}
In 2020, the COVID-19 pandemic transformed radically our way of living, communicating, socializing, commuting, shopping, as well as educating (\cite{Ali21,Camila21,Kostas21}). Most countries aimed to limit the spread of SARS-CoV-2 by non-pharmaceutical interventions (NPIs), including the implementation of lockdowns of varying intensity and geographic scope (\cite{ferguson2020}). Among so many lessons learned in these almost two years with COVID-19 (\cite{alamo2021}),  we know that vaccines are the major weapon to decrease the severity of COVID-19 (\cite{polack2020,dagan2021,henry2021}). Nevertheless, new cycles of infections are still taking place, saturating public health capacities (\cite{contreras2021}).

On the basis of the mathematical models, a crucial aspect for policymakers during the COVID-19 pandemic was to design public health policies to set up and remove social distancing measures. From a theoretical control perspective, this problem clearly falls into the classic framework of optimal control (\cite{lewis2012optimal,rawlings2017model}). To this end, simple objective functions and some constraints were defined to obtain a temporal sequence of measures that minimize the epidemic effects.
However, we have learned from previous works (\cite{sadeghi2020universal,federico2020taming,morris2021optimal,bliman2021best,ketcheson2020optimal,kohler2020robust}), that optimization problems arising from the control of epidemics by means of non-pharmaceutical measures are in essence complex and counter-intuitive, and the way to define objective functions and constraints is far from being trivial.

Optimal control deals with complex problems in which the major aim is to obtain the best possible performance by using the less possible control action (\cite{sethi2000,luenberger1979}), in a kind of cause-effect balance (\cite{bertsekas2012}).
However, in many applications such a cause-effect separation is not so clear. This is the case of SIR-type systems describing epidemics (models based on the seminal work \cite{kermack1927}), when the general objective is to minimize its effect/severity on the population. While it is clear that the number of infected individuals at any time and/or the total number of them at the end of the epidemic must be reduced, it is not so clear how to properly achieve such a goal. Evidence can be collected showing that the use of direct and simple control objectives - i.e., to minimize the number of infected individuals at all times or to minimize the final total number of infected - are not an option, since they may produce largely sub-optimal solutions. 

The key point to properly pose an optimal control problem (i.e., to define the objective function and the limits for the variables) is to understand its dynamic (including a stability analysis of the equilibrium and invariant regions) to know what can (and more importantly, what cannot) be done with the system under control. Asking for an impossible task, ignoring the limits that the system imposes, usually leads to poor performances. In the particular case of SARS-CoV-2, thinking that it is possible to fully control and/or manipulate a global epidemic by external intervention is, at best, a naive illusion and, at worst, a dangerous belief, since it motivates to take severe measures that produce the exact opposite effect to the one sought (as many examples around the world demonstrate, \cite{rypdal2020,abbasi2020,hale2021}). 

Several approaches have been proposed to find an optimal-control-based social distancing for SIR models in the context of the recent epidemic of COVID-19. The two main metrics to measure the disease impact are (\cite{di2021optimal}): the infected peak prevalence, IPP (maximal fraction of infected individuals along time), which is closely related to the health systems capacity, and the epidemic final size, EFS (total fraction infected).
A first result - when minimizing either EFS or IPP - is that optimal solutions can be obtained for the ubiquitous single-interval social distancing (\cite{sadeghi2020universal}), i.e., a fixed reduction of the reproduction number $\R$ (by reducing the infection rate) for a given period of time $\tau_f \!\!- \!\! \tau_s$, with $0 \! \leq \! \tau_s \! < \! \tau_f \! < \! \infty$.

In \cite{sadeghi2020universal,federico2020taming,morris2021optimal} rigorous analyses are made to show how to find the optimal single interval control action that minimizes IPP. The focus is put on the time at which optimal social distancing should start and finish and, as stated in \cite{morris2021optimal}, even when theoretical optimal (and near-optimal) social distancing are found, they are not robust to implementation errors, since "small errors in timing the intervention produce large increases in peak prevalence". The main problem of this strategy is, however, that it does not account for the other severity index, the EFS. So, the total number of infected (and the total number of deceases) is not minimized.

Similar approaches (and results) concerning single interval interventions are presented in \cite{bliman2021best,di2021optimal,ketcheson2020optimal}, but minimizing the EFS. In these works, it is found that optimal social distancing also exists, but in a rather unimplementable context. In \cite{bliman2021best} it is said that "the best policy consists in applying the maximal allowed social distancing effort until the end of the time interval, starting at a given particular time - that may be found by a simple algorithm". Indeed, the strategy consists in leaving the system in open-loop until the susceptible/non-infected fraction of individuals approaches the herd immunity threshold and then, at a particular time, implement the hardest possible social distancing. As it can be inferred, this strategy has two main drawbacks: any small error in the timing produces a performance drastically different from the optimal one and, more importantly, the infected peak prevalence is unacceptably large (since the system is left in open-loop for a long time before acting). The latter point is clearly mentioned in \cite{di2021optimal}, where it is said that "to minimize the total number of infected, the intervention should start close to the peak so that there is minimal rebound once the intervention is stopped". The former point, on the other hand, is demonstrated by means of simulation results in \cite{ketcheson2020optimal}, where slightly different social distancing from the optimal (or near-optimal social distancing) produce severe sub-optimal results. 

Another (practical) optimal control approach that closed the loop and led to more complex and realistic scheduling for social distancing can also be found in the recent literature. In works as \cite{kohler2020robust}, where it is proposed a model predictive control (MPC) based on the SIDARTHE model introduced in \cite{giordano2020sidarthe}, or \cite{morato2020optimal}, where the MPC is based on SIRASD (Susceptible-Infected-Recovered-Asymptomatic-Symptomatic-Dead) models, the control objective consists in minimizing both, the current number of infected individuals (and/or fatalities) and the time of isolation.  Other similar approaches concerning MPC strategies can also be found in \cite{alleman2020covid,peni2020nonlinear}.

In any case, the common factor in all the cited literature - at the best of the authors knowledge - is that no conclusive results are shown concerning which is the best policy to simultaneously minimize the IPP, EFS, and time social distancing lasts. It seems that always a trade off arises, which impeded the policymaker to reach a clear and simple consensus about the right measure for any possible scenario. In this article we show that the key point to achieve - or, at least, to arbitrarily approximate - such a goal is the way the optimization problem is posed, by properly separating transient and stationary regimes. More precisely, the objective of this work is to present - based on a pure dynamical analysis of the SIR-type models - a different perspective to formulate (pose) the social distancing optimal control problem. Instead of considering the control objective of minimizing the infected peak prevalence (IPP) or the epidemic final size (EFS), we just steer the susceptible to the (open-loop) herd immunity, since this threshold represents the minimal EFS at steady-state, for any finite-time intervention. Furthermore, by taking advantage of the fact that the infected peak is independent of the EFS, the IPP is maintained under an upper bound (computed in accordance with the health system capacity), while the only quantity to be minimized is the strength and time of the social distancing. As demonstrated by several simulation results, this strategy seems to be general enough to provide a confidence baseline to policymakers in the critical task of decision making in a pandemic context. 

\section{Review of control SIR Model}\label{sec:revSIR}
In this section we review the SIR epidemic model \cite{kermack1927}, which describes the fractions of susceptible $S(t)$ and infectious $I(t)$ individuals in a population at time $t$. New infections occur proportional to $S(t)I(t)$ at a transmission rate $\beta$, and infectious individuals recover or die at a rate $\gamma$. 

We consider also non-pharmaceutical interventions that reduce the effective transmission rate, $\beta(t)$, below its value in the absence of intervention, which is considered fixed. 
By rescaling the time by $\tau :=t \gamma$, the SIR model can be written in non-dimensional form as (\cite{sontag2011lecture,bertozzi2020challenges}):
\begin{subequations}\label{eq:SIRnondim}
\begin{align}
    \dot{S}(\tau) &= - \R(\tau) S(\tau) I(\tau), \\
    \dot{I}(\tau) &=   \R(\tau) S(\tau) I(\tau) - I(\tau),
\end{align}
\end{subequations}
where $\R(\cdot):=\beta(\cdot)/\gamma$ denotes the time-varying reproduction number fulfilling $\R (\cdot)\in \Omega_{\R}$, with
\begin{eqnarray}\label{eq:Rform}
	\Omega_{\R}\!\!:=\!\! \{\R(\cdot)\!:\! \mathbb R_{\geq 0} \rightarrow \mathbb R_{\geq 0} \!\!:\!\! \R(\tau) \!\!\in \!\![\uR,\oR], \mbox{~for~} \!\tau \!\!\in \!\![\tau_s,\!\tau_f], \nonumber\\
	\mbox{and~} \R(\tau) \!\!=\!\! \oR, \mbox{~for~} \!\tau \!\!\in\!\! [0,\!\tau_s) \mbox{~and~} \!\tau \!\!\in\!\! (\tau_f,\!\infty) \}, 
\end{eqnarray}
where $0\!<\!\tau_s\!<\! \tau_f \!<\! \infty$ denotes the starting and ending intervention time ($\tau_f$ is assumed to be finite to model the fact that social intervention has always an end), and $0\!<\!\uR\!<\! \oR$ are the minimal and maximal values for the reproduction number, i.e., $\oR$ is the value of the reproduction number in the absence of intervention and $\uR$ is the minimal value of the reproduction number corresponding to the maximal effectiveness of an intervention (the case $\uR\!=\!0$ is not considered, since a perfect full lockdown is impossible to implement). Later on, further conditions on the form of $\R(\cdot)$ will be established, to account for realistic interventions (as social distancing, face mask-wearing, etc.), that cannot be considered as a continuous function of time.

Susceptible $S$ and infectious $I$ are positive, and constrained to be in the set
\begin{eqnarray}
\setX \!:=\! \{(S,I) \in \mathbb R^2\!: S \in [0,1], I \in [0,1],S+I\leq 1 \}, \nonumber
\end{eqnarray}
for all $\tau \geq 0$. Particularly, denoting $\tau=0$ the epidemic outbreak time, it is assumed that $(S(0),I(0)):=(1-\epsilon,\epsilon)$, with $0 < \epsilon \ll 1$; \textit{i.e.}, the fraction of susceptible individuals is smaller than, but close to $1$, and the fraction of infected is close to zero at $\tau=0$.
\begin{rem}
	Note that even when compartments ($S$, $I$ and the removed, $R$) can be divided into sub-compartments, connected to each others to have a more detailed description of an epidemic (as in \cite{giordano2020sidarthe}) the main dynamic of the original system is maintained, as detailed in \cite{sadeghi2020universal}.
\end{rem}	
%
\subsection{No-intervention dynamical analysis}
%
We will assume first that $\R(\tau)\equiv \oR$, for $\tau \in [0,\infty]$, which represents the no-intervention (or open-loop) scenario. The solution of \eqref{eq:SIRnondim} for $\tau \geq \tau_0 >0$ - which was analytically determined in \cite{harko2014exact} - depends on $\oR$ and the initial conditions $(S(\tau_0),I(\tau_0)) \in \setX$. Since $S(\tau) \geq 0$, $I(\tau) \geq 0$, for $\tau \geq \tau_0 >0$, then $S(\tau)$ is a decreasing function of $\tau$ (by \ref{eq:SIRnondim}.a), for all $\tau \geq \tau_0$.
From \ref{eq:SIRnondim}.b, it follows that if $S(\tau_0)\oR \leq 1$, $\dot{I}(\tau) = (\oR S(\tau) -1)I(\tau) \leq 0$ at $\tau_0$. Furthermore, given that $S(\tau)$ is decreasing, $I(\tau)$ is also decreasing for all $\tau \geq \tau_0$.
On the other hand, if $S(\tau_0)\oR>1$, $I(\tau)$ initially increases, then reaches a global maximum, and finally decreases to zero. In this latter case, the peak of $I$ or infected peak prevalence IPP, $\hat I$, is reached at $\hat \tau$, when $\dot I = \oR SI - I=0$. $\hat I$ depends on initial conditions $S(\tau_0),I(\tau_0)$, and $\oR$, as follows:
\begin{eqnarray}\label{eq:Ipeak}
	IPP \!&:=&\! \hat I(\!\oR,\!S(\tau_0),\!I(\tau_0)) \nonumber \\
	&=&\!\! I(\tau_0) \!\!+\!\! S(\tau_0) \!\!-\!\! (1/\oR) (1\!\!+\!\!\ln(S(\tau_0) \oR)).
\end{eqnarray}
Condition $\dot I = \oR SI - I=0$ implies that $S=S^*$, where 
\begin{eqnarray}\label{eq:Sstar}
S^*:= \min \{1,1/\oR\} 
\end{eqnarray}
is a threshold or critical value, known as "herd immunity" (i.e., the value of $S$ under which $I$ cannot longer increase). This way, conditions $S(\tau_0)\oR>1$ and $S(\tau_0)\oR<1$ that determines if $I(\tau)$ increases or decreases at $\tau_0$ can be rewritten as $S(\tau_0)>S^*$ and $S(\tau_0)<S^*$, respectively.

For the sake of simplicity, we define $S_\infty:= \lim_{\tau\rightarrow \infty} S(\tau)$ and $I_\infty:= \lim_{\tau \rightarrow \infty} I(\tau)$, which are values that depend on initial conditions $S(\tau_0),I(\tau_0)$, and $\oR$. By taking $\tau \rightarrow \infty$ for the solutions proposed in \cite{harko2014exact}, we obtain $I_\infty = 0$. Furthermore, following a similar procedure than the one proposed in \cite{abuin2020characterization} for in-host models, $S_\infty$ fulfills the condition
\begin{eqnarray}\label{eq:Ssol1}
	-\oR S(\tau_0)e^{-\oR (S(\tau_0)+I(\tau_0))}= -\oR S_\infty  e^{-\oR S_\infty}.
\end{eqnarray}
Denote $y\!:=\!-\oR S_\infty$ and $z\!:=\!-\oR S(\tau_0)e^{-\oR (S(\tau_0)+I(\tau_0))}$. Then, equation \eqref{eq:Ssol1} can be written as $z=ye^y$, and $y$ can be obtained by $W(z)\!=\!y$, where $W(\cdot)$ is the Lambert function (\cite{pakes2015lambert}). That is, $W(-\oR S(\tau_0)e^{\!-\!\oR (S(\tau_0)+I(\tau_0))}) \!=\! -\oR S_\infty$, or
\begin{eqnarray}\label{eq:Sinf}
	S_\infty(\oR,S(\tau_0),\!I(\tau_0))\!\! :=\!\! -\frac{W(\!-\!\oR\! S(\tau_0)\!e^{-\oR (S(\tau_0)\!+\!I(\tau_0))})}{\oR}.\!\!
\end{eqnarray}
The epidemic final size, defined as $EFS \!:=\! 1- S_\infty$, is finally given by
\begin{eqnarray}\label{eq:Cinf}
	EFS = \!1 + \! (1/\oR)W(\!-\oR S(\tau_0) e^{-\oR (S(\tau_0)\!+\!I(\tau_0))}),
\end{eqnarray}
and, as the IPP, is a function of the initial conditions and $\oR$. 

The following Lemma states the maximum of function $S_\infty(\oR,S(\tau_0),I(\tau_0))$ over $\setX$
\begin{lem}[Maximum $S_\infty$ over $\setX$]\label{lem:Sinf_opt}
	Consider system \eqref{eq:SIRnondim} with arbitrary initial conditions $(S(\tau_0),I(\tau_0)) \in \setX$, for some $\tau_0 \geq 0$, and $\oR>0$ is fixed.
	Then, the maximum of $S_\infty(\oR,S(\tau_0),I(\tau_0))$ occurs at $(S^*,0)$ and it is given by $S^*$. Furthermore, if $I(\tau_0) \in [\delta,1]$, for some $\delta>0$, the maximum takes place at $(S^*,\delta)$ and it is given by $- W(- \oR S^* e^{-\oR (S^*+\delta)})/\oR$.
\end{lem}
\begin{pf} See Appendix 2.
\end{pf}

Next, some further properties regarding the general behaviour of $S_\infty(\oR,S(\tau_0),I(\tau_0))$, for different values of its arguments, are given. 
\begin{propt} \label{propt:sinfty}
	Consider system \eqref{eq:SIRnondim} with arbitrary initial conditions $(S(\tau_0),I(\tau_0)) \in \setX$, for some $\tau_0 \geq 0$, and $\oR>0$. Then:
	(i) $\lim_{\oR \rightarrow \infty} S_\infty(\oR,S(\tau_0),I(\tau_0))=0$ and $\lim_{\oR \rightarrow 0} S_\infty(\oR,S(\tau_0),I(\tau_0)) \approx S(\tau_0)$.
	(ii) For $S(\tau_0) > S^*$ and fixed $I(\tau_0)>0$ and $\oR>1$, $S_\infty(\oR,S(\tau_0),I(\tau_0))$ decreases with $S(\tau_0)$, and $S_\infty(\oR,S(\tau_0),I(\tau_0)) < S^*$. This means that the closer $S(\tau_0)$ is to $S^*$ from above, the closer will be $S_\infty$ to $S^*$ from below.
	(iii) For $S(\tau_0) < S^*$ and fixed $I(\tau_0)>0$ and $\oR>1$, $S_\infty(\oR,S(\tau_0),I(\tau_0))$ increases with $S(\tau_0)$, and $S_\infty(\oR,S(\tau_0),I(\tau_0)) < S^*$. This means that the closer $S(\tau_0)$ is to $S^*$ from below, the closer will be $S_\infty$ to $S^*$, from below. 
	(iv) For any fixed $S(\tau_0)$ and $\oR$, $S_\infty(\oR,S(\tau_0),I(\tau_0))$ decreases with $I(\tau_0)$ and $S_\infty(\oR,S(\tau_0),I(\tau_0)) \leq S^*$.
	(v) $\lim_{(S(\tau_0),I(\tau_0)) \rightarrow (S^*, 0)} S_\infty(\oR,S(\tau_0),I(\tau_0))=S^*$. 
	If $S(\tau_0) = S^*$ and $I(\tau_0) \approx 0$, $S_\infty \approx S^*$, for any value of $\oR$ (note that $S^*=1$ for $\oR<1$).
\end{propt}

\begin{figure}
	\centering
	\includegraphics[width=0.6\columnwidth]{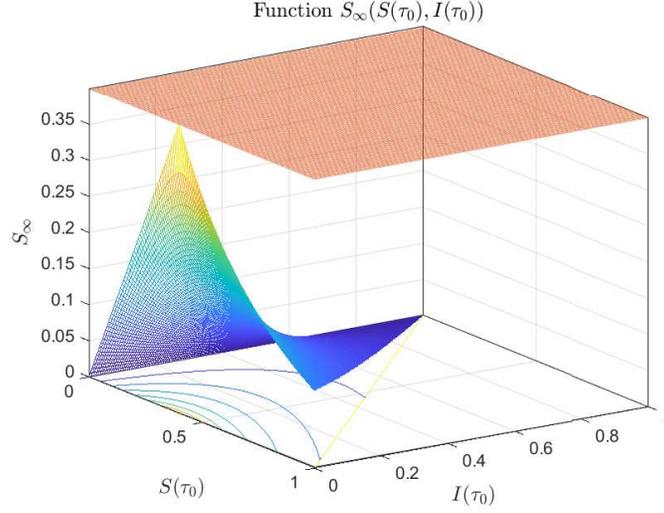}
	\caption{\small{Function $S_\infty(\R,S(\tau_0),I(\tau_0))$ is bounded from above by $S^*\!=\!1/\R$ ($S_\infty = S^*$, light red plane). Furthermore, $S_\infty$ reaches its maximum, given by $S^*$, at $S(\tau_0)\!=\!S^*$, $I(\tau_0)\!=\!0$.}}
	\label{fig:SinfFunc}
\end{figure}

The proof of Property \ref{propt:sinfty} is omitted for brevity. Figures \ref{fig:SinfFunc} and \ref{fig:Sinf_S} show how $S_\infty$ behaves for different initial conditions. 
%
%
\subsection{Equilibrium characterization and stability} \label{sec:eq_estabil}
%
The equilibrium of system \eqref{eq:SIRnondim}, with $\R(\tau)\equiv \oR$ for $\tau \in [0,\infty]$, is obtained by zeroing each of the differential equations. This way, for initial conditions $(S(\tau_0),I(\tau_0))\in \setX$, the equilibrium set is given by: 
\begin{eqnarray}
\setX_s \!\!:=\!\!\{(\bar S,\bar I) \in \setX\!\!:\! \bar S \in [0,S(\tau_0)], \bar I\!\!=\!\!0\}, \nonumber
\end{eqnarray}
since $S$ is decreasing.
Next, a key theorem concerning the asymptotic stability of a subset of $\setX_s$ is introduced.
\begin{theo}[Asymptotic Stability]\label{theo:stability}
Consider system \eqref{eq:SIRnondim} with $\R(\tau)\equiv \oR$ and constrained by $\setX$. Then, the set
\begin{eqnarray}
\setX_s^{st} \!\!:=\!\! \{(\bar S,\bar I) \in \setX\!\!:\! \bar S \in [0,S^*], \bar I\!\!=\!\!0\}, \nonumber
\end{eqnarray}
where $S^*$ is the herd immunity defined as $S^*:=\min \{1,1/\oR\}$, is the unique asymptotically stable (AS) set of system \eqref{eq:SIRnondim}, with a domain of attraction (DOA) given by $\setX$.  
\end{theo}
\begin{pf} See Appendix 2.
\end{pf}

A corollary of Theorem \ref{theo:stability}, concerning the properties of $\setX_s^{st}$, is presented next.

\begin{cor}\label{cor:Xst}
Consider system \eqref{eq:SIRnondim} with arbitrary initial conditions $(S(\tau_0),I(\tau_0))\in \setX$, for some $\tau_0 \geq 0$. Then: 
(i) Set $\setX_s^{st}$ is a subset of $\setX_s$ (for $\oR<1$, $\setX_s^{st} \equiv \setX_s$), and its size depends on $\oR$, but not on the initial conditions.
(ii) As detailed in Remark \ref{rem:nonas}, in Appendix 2, subsets of $\setX_s^{st}$ are $\epsilon-\delta$ stable but not attractive (\textit{i.e.}, even when $\setX_s^{st}$ is AS as a whole, no subset of it is AS). 
This is particularly true for the state $(S^*,0)$ which belongs to $\setX_s^{st}$.
(iii) If $\oR <1$, $S^* = 1$. Then, $\setX_s^{st} \equiv \setX_s$ and the so called healthy equilibrium $x_h:=(\bar S,0)$ with $\bar S = 1$, lies in $\setX_s^{st}$, and so it is $\epsilon-\delta$ stable, but not attractive (any small value of $I$ will make the system to converge to $(\bar S,0)$) with $\bar S < 1$. 
(iv) If $\oR > 1$, set $\setX_s$ can be divided into two sets, $\setX_s = \setX_s^{st} \cup \setX_s^{un}$, where 
$\setX_s^{un}:=\{(\bar S,\bar I) \in \setX: \bar S \in (S^*,1], \bar I=0\}$ is \textbf{an unstable equilibrium set} (which contains the healthy equilibrium). 
(v) Given that any compact set including an AS equilibrium set is AS, $\setX_s$ is AS, for any value of $\oR$. However, if $\oR>1$, it contains an unstable equilibrium set, $\setX_s^{un}$.
\end{cor}

\begin{figure}
	\centering
	\includegraphics[width=0.6\columnwidth]{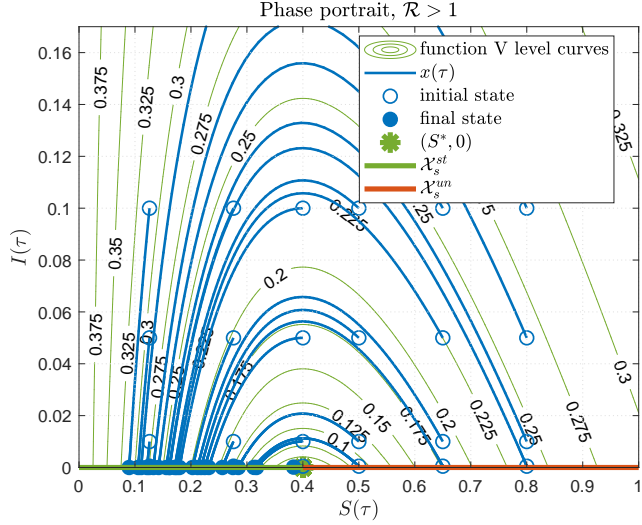}
	\caption{\small{Phase Portrait for system \eqref{eq:SIRnondim} with $\oR=2.9$ and initial conditions $(S(\tau_0),I(\tau_0)) \in \setX$. Set $\setX_s^{st}$ is in green (thick line), while $\setX_s^{un}$ is in red (thick line). Furthermore, the level curves of a Lyapunov function $V$ (similar to the one in \eqref{ec:lya1}, with $\bar x = (S^*,0)$) are plotted, which constitute invariant trajectories for the system (green thin line). As it can be seen, all the trajectories converge to $\setX_s^{st}$. Note that for $\oR<1$, $\setX_s^{st}$ is the whole set $\setX_s$.}}
	\label{fig:PhaPorRg1}
\end{figure}

Figures \ref{fig:PhaPorRg1} shows a Phase Portrait for system \eqref{eq:SIRnondim}, with $\oR\!>\!1$, and $(S(\tau_0),I(\tau_0)) \in \setX$. Similar behavior can be seen if $\oR<1$, with $\setX_s^{st} \!\equiv\! \setX_s$.
%
\section{Control}\label{sec:control}
Social distancing (or, more generally, non-pharmaceutical interventions) are the typical measures that policymakers implement to control epidemics when vaccination (effectiveness and distribution) is not enough. Social distancing lessens the disease transmission rate $\beta(\tau)$ or, directly, parameter $\R(\tau)$ in system \eqref{eq:SIRnondim} which, in turn, tends to reduce the two main indexes of the epidemic severity (\cite{di2021optimal}): the IPP and the EFS.

We assume now that $\R(\tau) \in \Omega_{\R}$. So, given that $\R$ is time-varying, IPP is no longer given by equation \eqref{eq:Ipeak}. However, given that the final intervention time $\tau_f$ is finite, equation \eqref{eq:Cinf} is still valid to describe EFS (i.e., $S_\infty$ is still governed by equation \ref{eq:Sinf}). The following Lemma gives an upper bound for the steady-state value of $S$ when $\R(\tau)$ belongs to $\Omega_{\R}$
%
%
\begin{lem}[$S$ steady-state upper bound] \label{lem:Suppbound}
	Consider system \eqref{eq:SIRnondim} with initial conditions $(S(0),I(0))=(1-\epsilon,\epsilon)$, $0 < \epsilon \ll 1$, and $\R(0)=\oR$ such that $S(0)\!>\!S^*$. Consider also that $\R(\cdot) \in \Omega_{\R}$. Then,
	(i) the system converges to an equilibrium state $(S_\infty,0)$ with $S_\infty=S_\infty(\uR,S(\tau_f),I(\tau_f))$ bounded from above by $S^*$, being $S^*<1$ the herd immunity corresponding to no social distancing and,
	(ii), the only way to achieve $S_\infty=S_\infty(\uR,S(\tau_f),I(\tau_f)) \approx S^*$ is with a $\R(\cdot) \in \Omega_{\R}$ producing $S(\tau_f)\approx S^*$ and $I(\tau_f)\approx 0$, which implies that system \eqref{eq:SIRnondim} achieves a Quasi steady-state (QSS) condition at $\tau_f$. Three particular cases where condition $S_\infty \approx S^*$ is not achieved are: 
	(ii.a) if $S(\tau_f) > S^*$ and $I(\tau_f)\approx 0$, then a \textbf{second outbreak wave} will necessarily take place at some time $\tau > \tau_f$ and, finally, the system will converge to $S_\infty< S^*$ (the greater is $S(\tau_f)$ with respect to $S^*$, the smaller will be $S_\infty$ with respect to $S^*$), 
	(ii.b) if $S(\tau_f) < S^*$ and $I(\tau_f)\approx 0$, then $S_\infty$ will be close to $S(\tau_f)$ (the smaller is $S(\tau_f)$ with respect to $S^*$, the smaller will be $S_\infty$ with respect to $S^*$), and
	(ii.c) if $I(\tau_f)$ does not approach zero (i.e., no QSS conditions is reached at $\tau_f$), then no matter which value $S(\tau_f)$ takes, $S_\infty$ will be smaller than $S^*$ (the farther is $S(\tau_f)$ from $S^*$ from above or from below, the smaller will be $S_\infty$ with respect to $S^*$).
\end{lem}
\begin{pf}
	See Appendix 2.
\end{pf}
%

	Lemma \ref{lem:Suppbound} is a simple but strong result concerning any kind of social distancing schedule (interrupted at some finite time). On one hand, Lemma \ref{lem:Suppbound} (i) establishes that the minimal possible EFS is completely determined by the epidemic itself (the value of $\oR$) and, provided that no immunization (by vaccination and/or development of the individual's immune system) is considered, it cannot be modified by non-pharmaceutical measures. Then, any proposed social distancing schedule must be able to achieve this minimal EFS since, otherwise, values of $S_\infty$ (significantly) smaller than $S^*$ will produce more total infected individuals (and fatalities) than possible.
	On the other hand, Lemma \ref{lem:Suppbound} (ii) establishes that $S^*$ must be reached as a QSS condition (i.e., with $I(\tau_f) \! \approx \!0$), since otherwise $S(\tau)$ will decrease after the social distancing is interrupted at $\tau_f$, so producing $S_\infty \!< \!S^*$.
%
\begin{exmp} To clearly show this later property, we selected a particular problem (introduced in \cite{bliman2021best}), consisting in system \eqref{eq:SIRnondim}, with $\oR=2.9$ ($\beta=0.29$ days$^{-1}$ and $\gamma=0.1$ days$^{-1}$), $I(0)=1.49 \times 10^{-5}$, $S(0) = 1-I(0)$ and $\uR=0.66$. These settings will be used throughout the paper to demonstrate the proposed control strategies. 
Figure \ref{fig:PhaPor1} shows the phase portraits of this system under different control strategies, $\R(\cdot) \in \Omega_\R$. In Figure \ref{fig:PhaPor1} (left), two controls $\R(\cdot) \in \Omega_\R$ are applied (one strong and the other soft) at some time $\tau_s$, being $S(\tau_s) = 0.9$ and $I(\tau_s) = 0.065$, up to a large enough time $\tau_f$, such that the system reaches a QSS. In the first case (blue line) $I(\tau_f)\approx 0$ and $S(\tau_f)=0.7$. Given that $S(\tau_f)$ is significantly greater than $S^*=0.3448$ and $I(\tau_f)$ is small but positive, a second wave occurs (which explains the dynamic behaviour experienced when hard social distancing are implemented for a long period of time), that steers $S_\infty$ to a value significantly smaller than $S^*$ ($S_\infty = 0.13$). In the second case (red line), $I(\tau_f)\approx 0$ and $S(\tau_f)=0.15$. Since $S(\tau_f)$ is significantly smaller than $S^*$, and it cannot longer grow for $\tau > \tau_f$, so $S_\infty$ is again significantly smaller than $S^*$ ($S_\infty = 0.15$). The key point is that, although $\tau_f$ is large, it is always finite, so the state evolution is governed by the open-loop behavior for $\tau > \tau_f$. This can be seen by considering the level curves of a Lyapunov function for the open-loop system (green lines), which encircles the equilibrium $(S^*,0)$. These level curves are invariant trajectories for the open-loop system, so once the system reaches one of them at $\tau_f$, it will continue on the same curve for all $\tau>\tau_f$, in such a way that outer states at $\tau_f$ correspond to outer states for $\tau \rightarrow \infty$.

In Figure \ref{fig:PhaPor1} (right), the same two controls are applied at the same time $\tau_s$, but they are interrupted at some finite time $\tau_f$ such that the system has not reached a QSS condition. In the first case (blue line) the state evolution is such that $S(\tau_f)=0.8$ and $I(\tau_f)=0.045$. Then, for $\tau>\tau_f$, the system goes along one level curve of the Lyapunov function, which steers the system to $(S_\infty,I_\infty) \approx (0.09, 0)$, with $S_\infty$ significantly smaller than $S^*$. In the second case (red line), $S(\tau_f)=0.49$ and $I(\tau_f)=0.2$. This way, the corresponding level curve of the Lyapunov function steers the system - again - to $(S_\infty,I_\infty) \approx (0.09, 0)$, with $S_\infty$ significantly smaller than $S^*$, showing that any control action interrupted before a QSS is reached will necessarily produce a steady-state with $S_\infty$ smaller than $S^*$.
\begin{figure*}
	\centering
	\begin{subfigure}{.4\textwidth}
		\centering
		\includegraphics[width=1 \columnwidth]{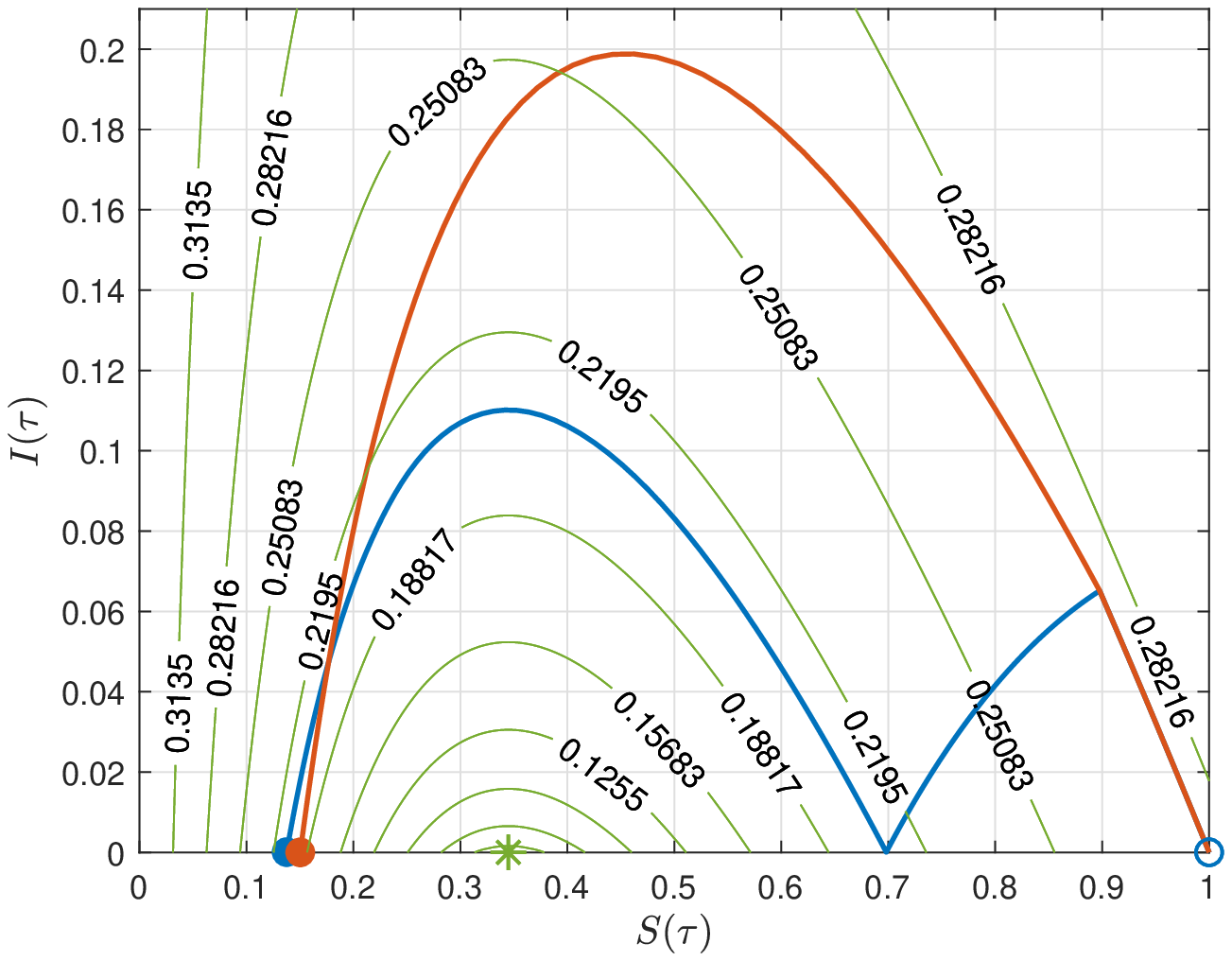}
		\label{fig:PhaPor1a}
	\end{subfigure}%
	\hspace*{0.2truecm}
	\begin{subfigure}{.4\textwidth}
		\centering
		\includegraphics[width=1 \columnwidth]{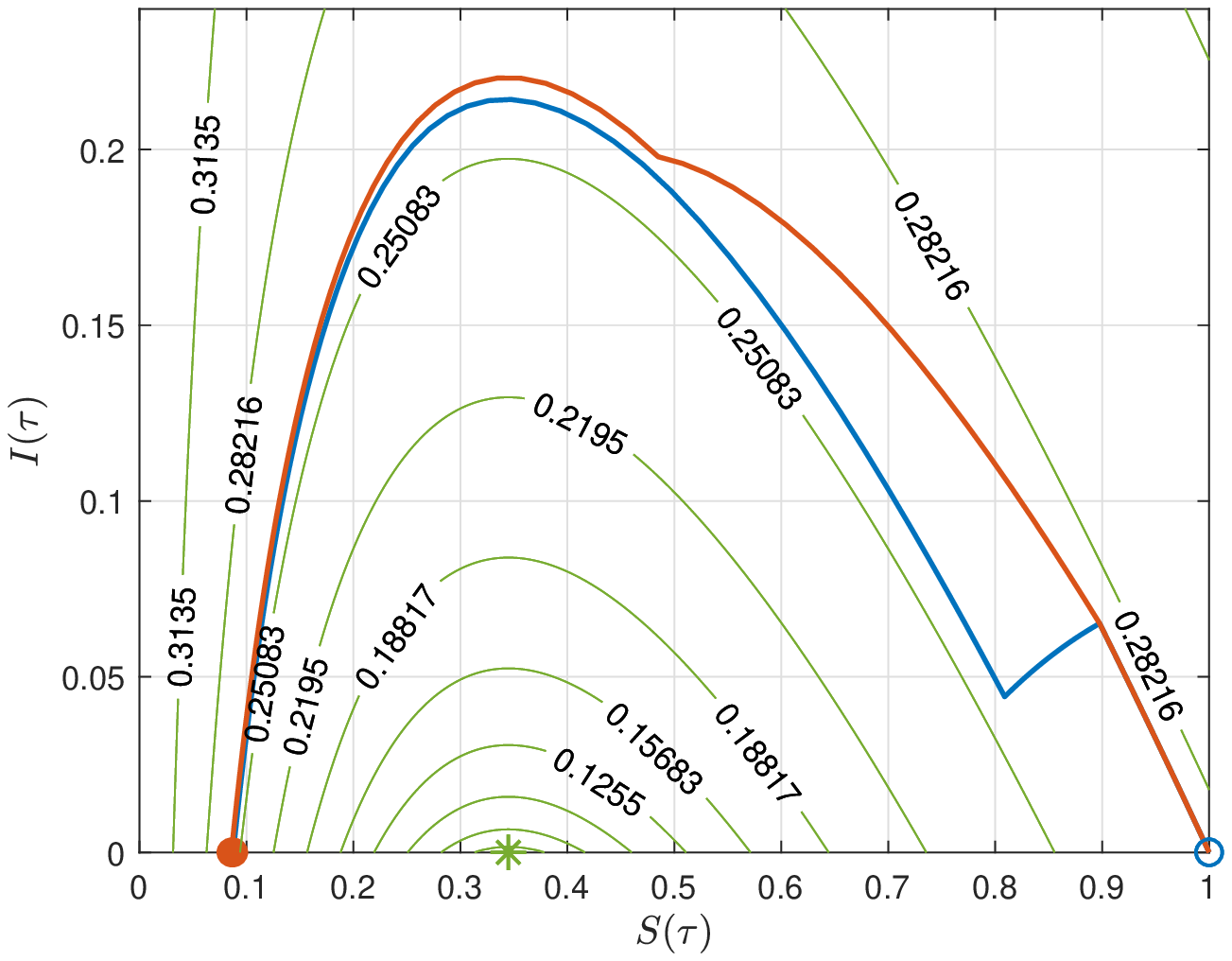}
		\label{fig:PhaPor1b}
	\end{subfigure}
	\caption{\small{Long term (left) and short term (right) control. No matter which control $\R(\cdot) \in \Omega_{\R}$ is applied, if it is interrupted at some finite time $\tau_f$ such that the system has not reached the quasi steady-state condition $(S(\tau_f),I(\tau_f))\! \approx \!(S^*,0)$, condition $S_\infty \! \approx \!S^*$ will not be achieved. Initial state: empty circle, $(S_\infty,\! I_\infty)$: filled circles, level curves of the Lyapunov function $V(S,I) \!:=\! S\!- \!S^* \!- \!S^* \ln(\frac{S}{S^*})\! + \!I$: green lines.}}
	\label{fig:PhaPor1}
\end{figure*}
\end{exmp}

Now, the question that naturally arises is whether or not it is possible to arbitrarily reduce the IPP, while maintaining $S_\infty \approx S^*$. To have a first insight into the answer, consider the integral of equation \eqref{eq:SIRnondim}.b, $I(\tau) = \int_{0}^{\tau} \R(t) S(t)I(t) dt  - \int_{0}^{\tau} I(t) dt + c$,
%
%
where $c$ is a constant determined by the initial values $S(0)$ and $I(0)$. Then, $\int_{0}^{\tau} I(t) dt \!\!= \!\!\int_{0}^{\tau} \R(t) S(t)I(t) dt  \!\!- \!\!I(\tau) \!\!= \!\!-S(\tau) \!\!-\!\! I(\tau)\!\! + \!\!c$,
%
%
and taking $\tau \!=\!0$, it follows that $c\!=\!S(0)\!+\!I(0)\!=\!1\!-\! \epsilon \!+\! \epsilon \!=\! 1$.
Now, taking the limits for $\tau \! \rightarrow \! \infty$, and recalling that $\R(\tau) \! \equiv \! \oR$ for $\tau \! \in \! (\tau_f,\infty)$, it follows that
\begin{eqnarray*}
	\int_{0}^{\infty} I(t) dt = 1 -S_\infty - I_\infty = 1 -S_\infty.
\end{eqnarray*}
This latter equality means that, even when $\R$ varies over time, $S_\infty$ only determines the area under the curve of $I(\tau)$, $AUC_I \!:=\! \int_{0}^{\infty} \!I(t) \!dt$, but not its peak $\hat I$. In other words, it is possible to minimize the EFS and also keep the IPP under a maximal value imposed by the health system capacity, as long as we respect the condition $\int_{0}^{\infty} \!I(t) dt \!=\! 1\! - \!S^*$.
%
%
\begin{rem}
	According to the latter discussion, it is apparent that it is not right for an optimal control strategy to simply minimize $\int_0^{\infty} I(t) dt$, as it is often done in the literature. Indeed, such a selection does not lead to a well-posed optimal control problem, since an unnecessary competition between steady-state objectives (minimize the EFS) and transient ones (minimize the IPP by minimizing $I$) arises, which necessarily produces an intermediate solution that optimizes neither EFS nor IPP. This fact can be clearly seen in \cite{bliman2021best,ketcheson2020optimal,di2021optimal} where the theoretical maximal value of $S_\infty=S^*$ is reached but the peak of $I$ is not controlled and, reciprocally, in \cite{federico2020taming,morris2021optimal,morato2020optimal,kohler2020robust}, where the IPP is controlled, but the EFS is always significantly greater than $1-S^*$. 
\end{rem}
In this work, we propose to pose the control objective in a rather different way. We divide the objectives into a primary or pure epidemiological one, and a secondary one, consisting in minimizing the side effects of the social distancing implemented to achieve the primary objective. 
The pure epidemiological control objectives read:
\begin{defn}[Epidemiological Control Objective]\label{defn:epidobj}
	Consider system \eqref{eq:SIRnondim} with initial conditions $(S(0),I(0))=(1-\epsilon,\epsilon)$, $0 < \epsilon \ll 1$, and $\R(0)=\oR$ such that $S(0)\!>\!S^*$. Consider also that $\R(\cdot) \in \Omega_{\R}$, and that a maximal value for $I$, $I_{max}>0$, is established according to the health system capacity. Then, the epidemiological control objective consists in steering $S(\tau)$ to $S^*$, as $\tau \rightarrow \infty$, while maintaining $I(\tau) \leq I_{max}$, for all $\tau \in \mathbb R_{\geq 0}$.
\end{defn}
The secondary objective accounts for the social distancing measures severity and, implicitly, for all the collateral effects of this intervention (economic, educational, and psychological problems, among others):
\begin{defn}[Social/Economic Control Objective]\label{defn:econobj}
	Consider system \eqref{eq:SIRnondim} with initial conditions $(S(0),I(0))=(1-\epsilon,\epsilon)$, $0 < \epsilon \ll 1$, and $\R(0)=\oR$ such that $S(0)\!>\!S^*$. Consider also that $\R(\cdot) \in \Omega_{\R}$, and that a maximal value for $I$, $I_{max}>0$, is established according to the health system capacity. Then, the social/economic control objective consists in minimizing $\int_{0}^{\infty} \oR - \R(t) dt$, provided that the Epidemiological Control Objective is achieved.
\end{defn}
To measure the social distancing severity, we define the Social Distancing Index, as follows:
\begin{eqnarray} \label{eq:sdi}
	SDI := \int_{0}^{\infty} \oR - \R(t) dt.
\end{eqnarray}
%
%
\section{Discussion}\label{sec:disc}
%
We now resume the fundamental question about the finding of some $\R(\cdot) \in \Omega_\R$ that fulfills the Epidemiological Control Objective \ref{defn:epidobj}. We proceed from the simplest to the most complex case, in which not only the epidemiological objective is considered but also the social/economic one.
%
\subsection{Single interval social distancing}\label{sec:sisd} 
%
Here we show that, under some mild conditions, it is possible to find a single interval social distancing that produces both $S_\infty\approx S^*$ and $I(\tau) \leq I_{max}$, for all $\tau \in \mathbb R_{\geq 0}$. A single interval social distancing is defined as
%
%
%
\begin{eqnarray} \label{eq:sisd}
	\R_{si}(\tau) := \left\{ 
	\begin{array}{cc}
		\oR  &  \mbox{for}~ \tau \in [0,\tau_s) \cup  (\tau_f,\infty) \\
		\R_{si}   & \mbox{for}~ \tau \in [\tau_s,\tau_f]
	\end{array} \right.
\end{eqnarray}
where $\R_{si} \in [\uR,\oR]$ is a fixed value of social distancing. Clearly, $\R_{si}(\tau) \in \Omega_\R$. This kind of control action was accounted for in \cite{bliman2021best,di2021optimal,ketcheson2020optimal} with the only objective of achieving $S_\infty=S^*$, and by \cite{morris2021optimal,sontag2021explicit}, with the only objective of minimizing $IPP$, respectively. 
%

\begin{figure}
	\centering
	\includegraphics[width=0.6 \columnwidth]{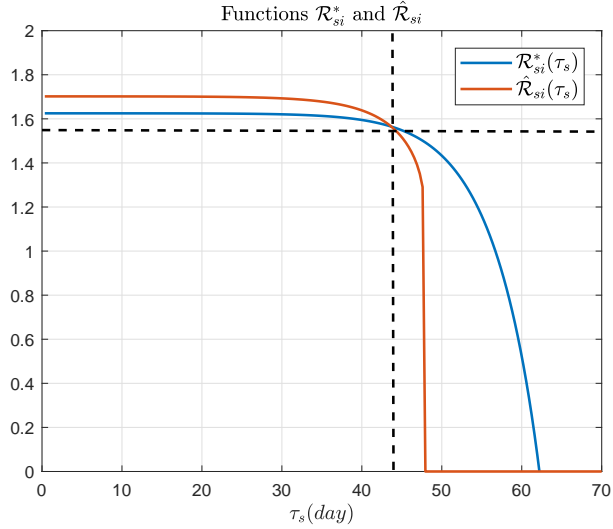}
	\caption{\small{Functions $\R_{si}^*(S^*,\tau_s)$ and $\hat{\R}_{si}(I_{max},\tau_s)$, for $\oR=2.9$ and $I_{max}=0.1$. Both reproduction numbers coincides at $\tau_s^g = 43.71$ days, where $\R_{si}^g = 1.565$}}
	\label{fig:Rs}
\end{figure}

In a more general context, it is possible to find a fixed reproduction number $\R^*_{si} \! \in \! [\uR,\oR]$ that, for a given social distancing starting time $\tau_s$, and a large enough $\tau_f$, produces $S_\infty \! \approx \!S^*$. By making $S_\infty(\R_{si}^*,\!S(\tau_s),\!I(\tau_s)) \!= \!S^*$, we have,
\begin{eqnarray} \label{eq:Rstar}
	\R_{si}^*(S^*,S(\tau_s),I(\tau_s)) := \frac{\ln S(\tau_s) - \ln S^*}{S(\tau_s)+I(\tau_s)-S^*}.
\end{eqnarray}
Given that previous to time $\tau_s$, system \eqref{eq:SIRnondim} evolves in open-loop, then $(S(\tau_s),I(\tau_s))$ are determined by the (unique) solution of system \eqref{eq:SIRnondim} under the effect of $\oR$. So we can write $\R_{si}^*(S^*,\tau_s)$ to emphasise the functional relationship between $\R_{si}^*$ and $\tau_s$. Indeed, for a given $S^*>1$, $\R_{si}^*(S^*,\tau_s)$ is a decreasing function of $\tau_s$, as it is shown in Figure \ref{fig:Rs}, blue line. Furthermore, $\R_{si}^*(S^*,0) \approx -\frac{\ln S^*}{1-S^*}$ (given that $S(0) \approx 1$ and $I(0) \approx 0$) and $\R_{si}^*(S^*,\hat \tau) =\frac{\ln S^* - \ln S^*}{S^*+\hat I-S^*}= 0$, since once the open-loop $S$ evolution reaches $S^*$ at $\hat \tau$, with a high value of $\hat I(\!\oR,\!S(\tau_0),\!I(\tau_0))=I(\hat \tau)$, nothing can be done to reach the condition $S_\infty \approx S^*$.

Now we ask for a fixed reproduction number, $\hat{\R}_{si}\! \in \![\uR,\oR]$, that guarantees $I(\tau) \! \leq \! I_{max}$, for all $\tau \! \in \! \mathbb R_{\geq 0}$. By making $\hat I(\hat{\R}_{si},S(\tau_s),I(\tau_s))\!=\!I_{max}$, we obtain the implicit function 
\begin{eqnarray} \label{eq:Rhat}
	\hat{\R}_{si}=\hat{\R}_{si}(I_{max},\tau_s).
\end{eqnarray}
For a given $I_{max}$, $\hat{\R}_{si}(I_{max},\tau_s)$ is a decreasing function of $\tau_s$, as it is shown in Figure \ref{fig:Rs}, red line.

Finally, by merging the latter two conditions, it is possible to define a (so-called) goldilocks social distancing:
\begin{defn}[Goldilocks single interval social distancing]
	The goldilocks single interval social distancing is defined by a starting time, $\tau_s^g$, fulfilling condition $\R_{si}^*(S^*,\tau^g_s) = \hat{\R}_{si}(I_{max},\tau^g_s)$,
	%
	%
	and the fixed reproduction number value, $\R_{si}^g := \R_{si}^*(S^*,\tau^g_s) = \hat{\R}_{si}(I_{max},\tau^g_s)$.
\end{defn}
The goldilocks single interval social distancing allows us to establish the following Theorem:
\begin{theo}\label{theo:gold}
	Consider system \eqref{eq:SIRnondim} with initial conditions $(S(0),I(0))=(1-\epsilon,\epsilon)$, $0 < \epsilon \ll 1$, and $\R(0)=\oR$ such that $S(0)\!>\!S^*$. Consider a given $I_{max}$ and consider also that $\R(\cdot) \in \Omega_{\R}$. Then, if for $S^*$ and $I_{max}$ there exists a goldilocks single interval social distancing, it is the only one that arbitrarily approaches the epidemiological control objectives, as $\tau_f \rightarrow \infty$. 
\end{theo}
\begin{pf}
	See Appendix 2.
\end{pf}
\begin{exmp}
Now we resume the example introduced in the later section, to show the system evolution when the goldilocks single interval social distancing is applied. We consider $I_{max}=0.1$, while $S^*=0.3448$. In Figure \ref{fig:Rs}, the social distancing starting time and reproduction number correspond to the goldilocks scenario ($\tau_s^g = 43.71$ days and $\R_{si}^g = 1.565$), respectively. Figure \ref{fig:GoldSI} (left), shows $S(\tau)$ (blue, upper plot), $I(\tau)$ (red, upper plot) and $\R(\tau)$ (lower plot) for a period of time of $300$ days. On the other hand, Figure \ref{fig:GoldSI} (right) shows the corresponding phase portrait, and the level curves of Lyapunov functions, $V^*$ (green curves) and $V^g$ (red curves), corresponding to $\oR$ and $\R_{si}^g$, respectively. The level curves of $V^g$ crosses that of $V^*$ and, given that $\R_{si}^g < \oR$, there is one that 'guides' the system exactly to $(S^*,0)$. This level curves is the one picked by the goldilocks single interval social distancing, applied at $\tau_s^g$. 
\begin{figure*}
	\centering
	\begin{subfigure}{.4\textwidth}
		\centering
		\includegraphics[width=1 \columnwidth]{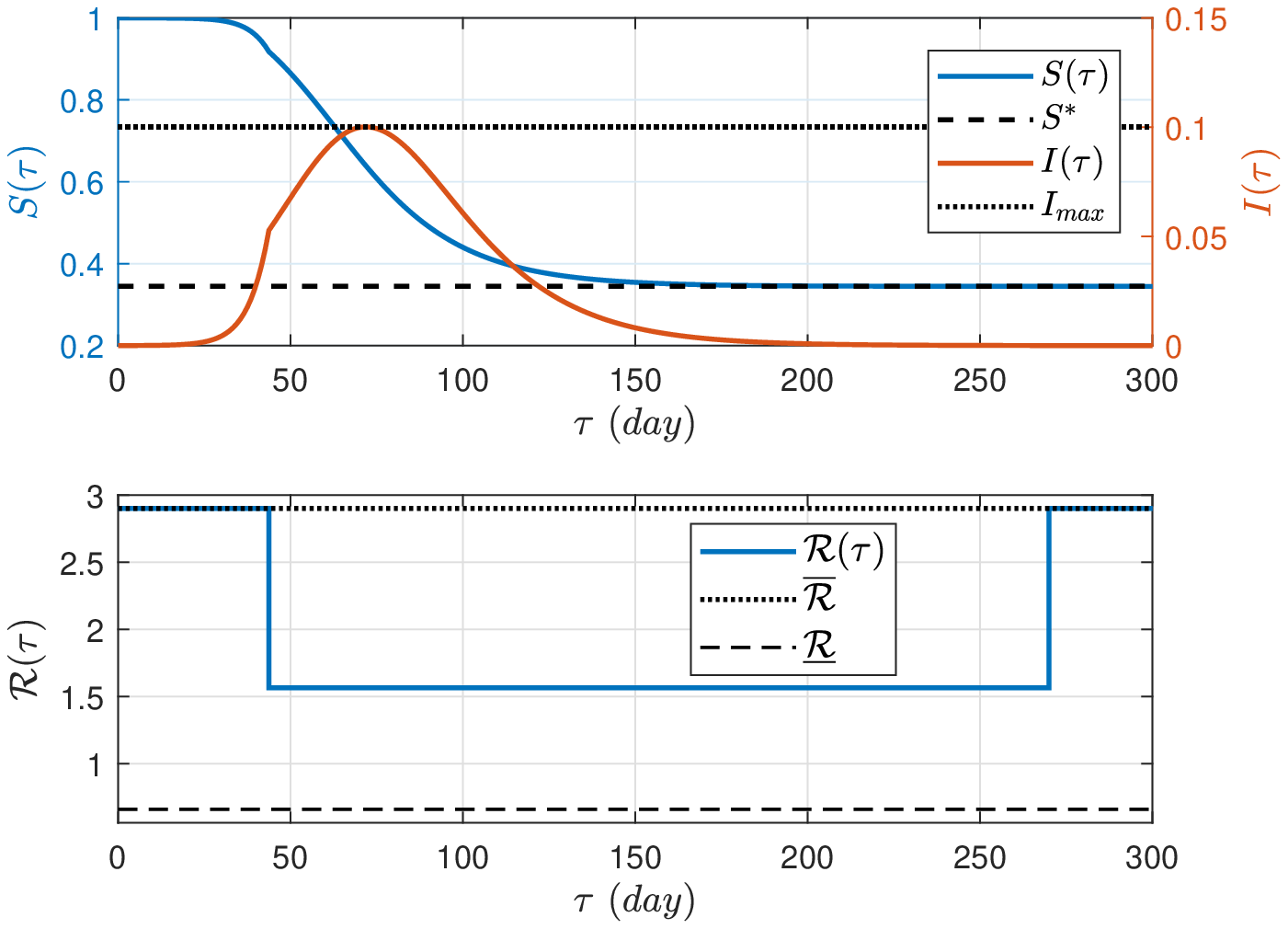}
	\end{subfigure}%
	\hspace*{0.2truecm}
	\begin{subfigure}{.4\textwidth}
		\centering
		\includegraphics[width=1 \columnwidth]{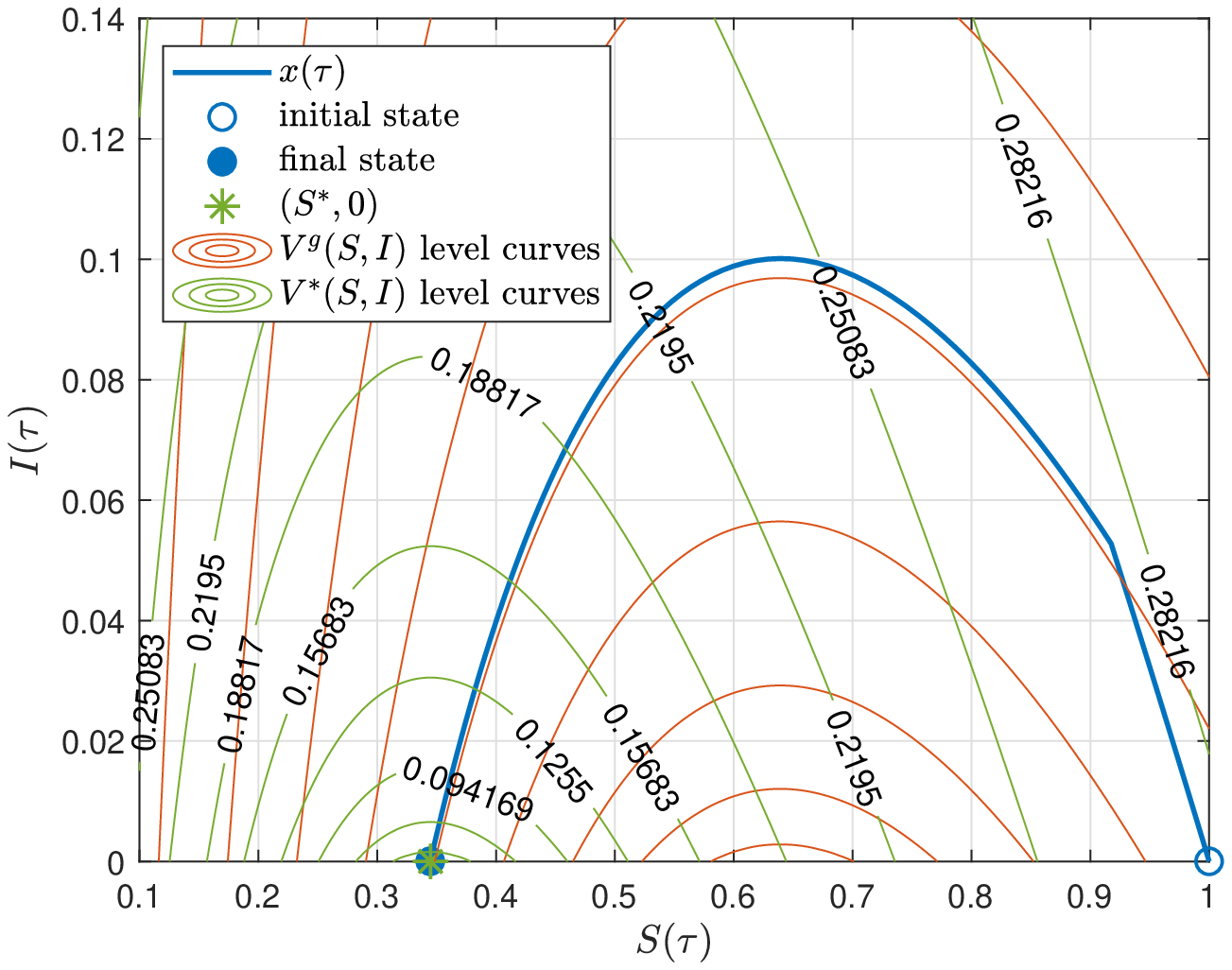}
	\end{subfigure}
	\caption{\small{$S,V$ and $\R$ time evolution, (left), and phase portrait, (right). System with $S^*=0.3448$ and $I_{max}=0.1$, when the goldilocks single interval social distancing is applied. The indexes are: EFS$=0.6600$, IPP$=0.1001$, SDI$=301.7100$.}}
	\label{fig:GoldSI}
\end{figure*}

The existence of the goldilocks single interval social distancing depends on the values of $S^*$ and $I_{max}$, and will not be analyzed here, for the sake of brevity. Anyway, goldilocks social distancing should be understood just as a theoretical approach, since it cannot be applied in realistic cases. On one side, $S_\infty$ and $\hat I$ are too sensitive to the values of $\tau^g_s$ and $\R_{si}^g$, so small changes in these values produce values of $S_\infty$ significantly smaller than $S^*$ and values of $\hat I$ significantly greater than $I_{max}$ (note that the exact value $\R^g_{si}$ may be difficult - if not impossible - to achieve/apply in a realistic setting intervention). On the other side, the goldilocks social distancing consists in a significantly low value of $\R^g_{si}$, applied during an extremely long time $\tau_f$ (producing severe side effects). The social distancing index corresponding to the goldilocks single interval social distancing is $SDI$=301.7100.
%

\end{exmp}
%
\subsection{'Wait, maintain, suspend' strategy} \label{sec:morcomplex} 
%
There is another (more complex) strategy to account for the epidemiological control objectives of Definition \ref{defn:epidobj}, that avoids the problem of the existence of a solution for any $I_{max}$. We will denote it - following the ideas presented in \cite{morris2021optimal} - 'wait, maintain, suspend' (the social distancing) strategy:
%
%
%
\begin{eqnarray} \label{eq:wms}
	\R(\tau) := \left\{ 
	\begin{array}{cc}
		\oR  &  \mbox{for}~ \tau \in [0,\tau_s) \cup (\tau_f,\infty),   \\
		\frac{1}{S(\tau)}   & \mbox{for}~ \tau \in [\tau_s,\tau_1),\\
		\R_{si}^*   & \mbox{for}~ \tau \in [\tau_1,\tau_f],
	\end{array} \right.
\end{eqnarray}
where $0 \! \leq \! \tau_s \! < \! \tau_1 \! < \! \tau_f \! \leq \! \infty$, $\tau_1$ is the time at which a threshold condition (specified later on) is reached, and $\R_{si}^* \! =\!  \R_{si}^*(S^*,\tau_1)$ is the fixed social distancing that, if started at $\tau_1$ and $\tau_f$ is large enough, produces $S_\infty \!  \approx S^*$ and $I_\infty \! =\! 0$ ($\tau_f$ large enough means the system reaches a QSS condition).

Time $\tau_s$ is considered now as the time at which the open-loop system reaches $I(\tau_s)\!=\!I_{max}$. At this time, the control action $\R(\tau)\!=\!\frac{1}{S(\tau)}$ is applied to system \eqref{eq:SIRnondim}, making $I(\tau)$ constant for the period $[\tau_s,\tau_1)$. As a result, $S(\tau)$ decreases linearly for $[\tau_s,\tau_1)$ (since $\dot S(\tau) \!= \!I(\tau) \!=\! I_{max}$). Now, if time $\tau_1$ is not large enough, $I(\tau)$ may increases for $\tau \! \geq \! \tau_1$, reaching a peak that overpasses $I_{max}$, violating this way the control objective $I(\tau) \! \leq \! I_{max}$. On the other hand, if $\tau_1$ is too large, $S(\tau)$ may decreases under $S^*$, violating the control objective $S_\infty \! \approx \! S^*$. This means that some $\tau_1$, denoted as $\tau_1^*$, exists such that both conditions are fulfilled (see the proof of Theorem \ref{theo:mcsol}, in Appendix 2), and the epidemiological control objective is arbitrarily approached.

The next Theorem, formalizes the benefits of the 'wait, maintain, suspend' strategy.
\begin{theo}\label{theo:mcsol}
	Consider system \eqref{eq:SIRnondim} with initial conditions $(S(0),I(0))\!=\!(1\!-\!\epsilon,\epsilon)$, $0 \! < \! \epsilon \! \ll \! 1$, and $\R(0)=\oR$ such that $S(0)\!>\!S^*$. Consider a given $I_{max}$ and consider also that $\R(\cdot) \in \Omega_{\R}$. Then, there exists some $0\! \leq \! \tau_s \! < \! \tau_1^* \! < \! \tau_f$ such that the "wait, maintain, suspend" strategy given by \eqref{eq:wms}, produces $(S(\tau),I(\tau))$ for $\tau \! \geq \! 0$ that arbitrarily approach the epidemiological control objectives, as $\tau_f \! \rightarrow \! \infty$. 
\end{theo}
\begin{pf}
	See Appendix 2.
\end{pf}
\begin{exmp}
Figure \ref{fig:WMS} shows the system evolution when the "wait, maintain, suspend" social distancing is applied to the system previously defined. As shown in Figure \ref{fig:Rwms}, in Appendix 2, $\tau_s = 47.8$ days, $\tau_1^* = 68.7$ days, and $\R_{si}^*(S^*,\tau_1^*) = 1.5654$. As it can be seen, the epidemiological control objectives are reached, since $S_\infty \approx S^*$ and $I(\tau)\leq I_{max}$ for $\tau \in [0,300]$ days. However, as before, the social distancing schedule is clearly unrealistic: the control action applied in the interval $[\tau_s,\tau_1^*]$ varies continuously in time (which impedes its implementation), and the severity of the social distancing is extremely high: the social distancing index is given by SDI$=298.8647$. Furthermore, similar to the goldilocks single interval social distancing, it can be shown that this strategy is extremely sensitive to wrong choice of $\tau_1^*$.
\begin{figure*}
	\centering
	\begin{subfigure}{.4\textwidth}
	\centering
	\includegraphics[width=1 \columnwidth]{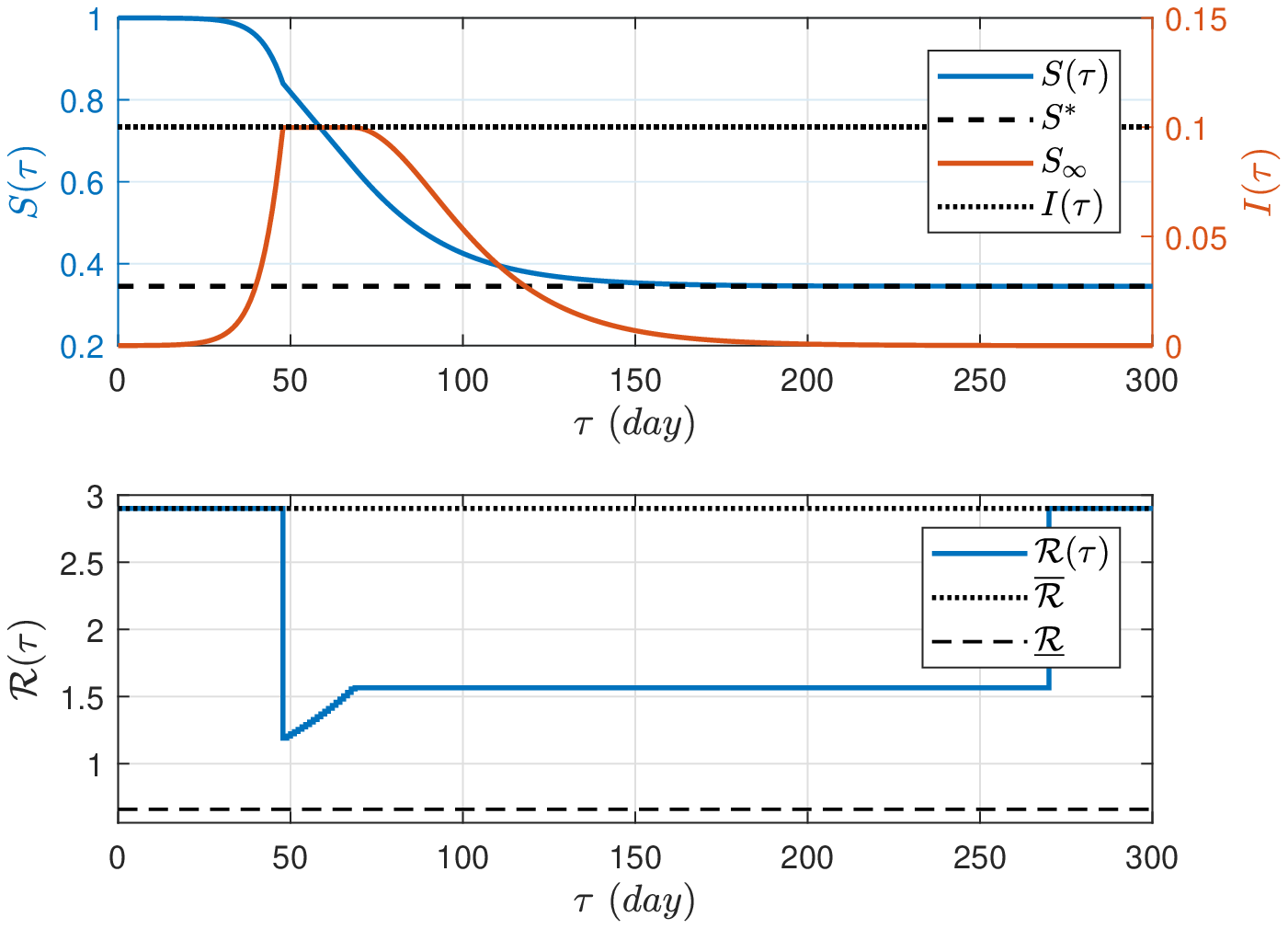}
\end{subfigure}%
\hspace*{0.2truecm}
\begin{subfigure}{.4\textwidth}
		\centering
	\includegraphics[width=1 \columnwidth]{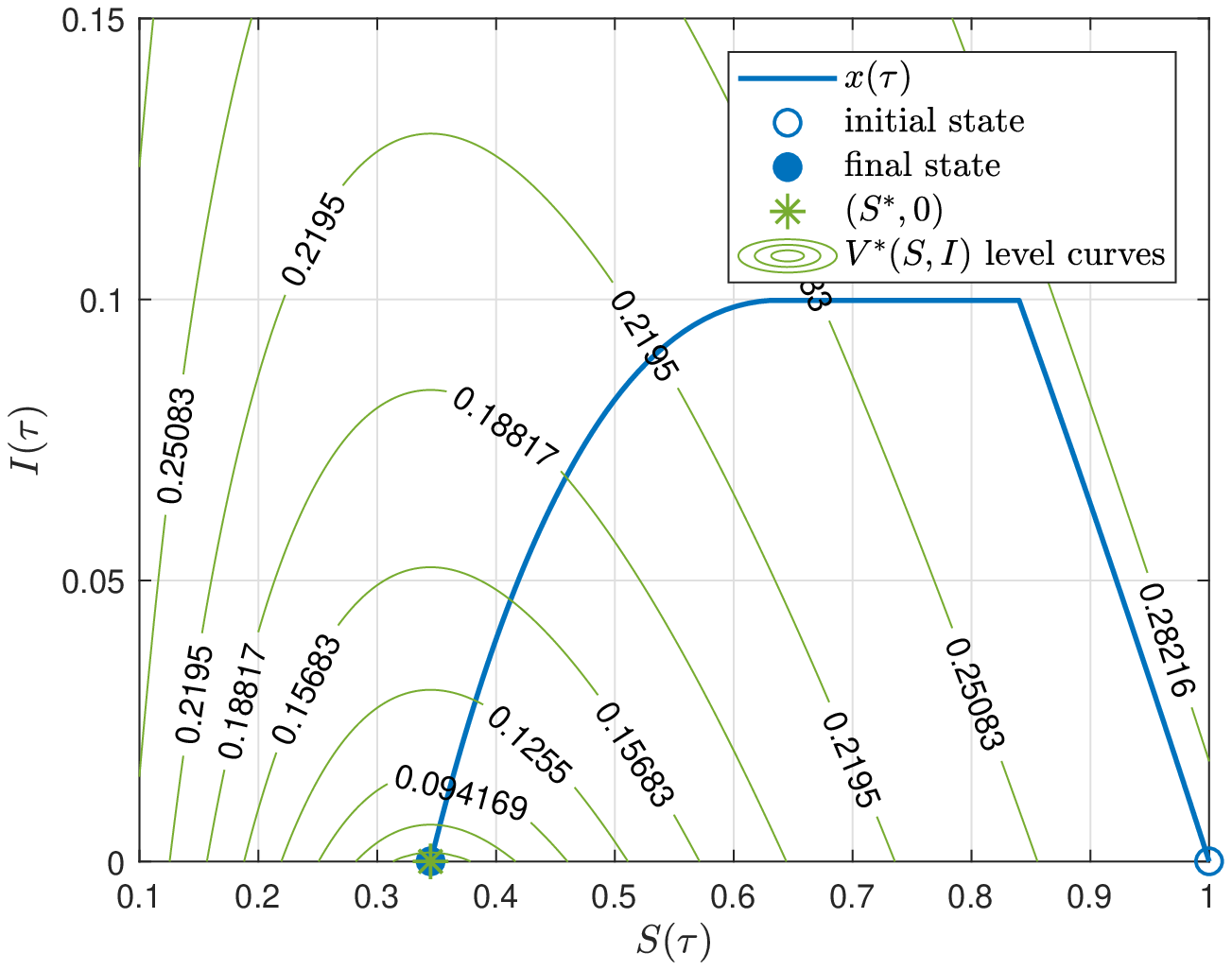}
		\end{subfigure}
\caption{\small{$S,V$ and $\R$ time evolution, (left), and phase portrait, (right). System with $S^*=0.3448$ and $I_{max}=0.1$, when the 'wait, maintain, suspend' control sequence \eqref{eq:wms} is applied. The corresponding performance indexes are: EFS$=0.6596$, IPP$=0.0998$, SDI$=298.8647$.}}
\label{fig:WMS}
\end{figure*}
%
%
\end{exmp}
%
\subsection{Optimal control strategy}\label{sec:optcont}  
%
Now, by taking advantage of the analysis of the previous strategies, and considering that any realistic control strategy should also consider the Social/Economic Control Objective of minimizing the social distancing severity, the following optimal control problem, $\mathcal P_{opt}(S(0),I(0),S^*,I_{max};\R(\cdot))$, is proposed:
\begin{eqnarray} 
	\begin{array}{rlr}
		\min\limits_{\R(\cdot)}\;\!\!
		&V(\R(\cdot)) = \int_{0}^{T} \oR -\R(t) dt & \\
		s.t.&&   \\
		&\dot S(\tau) = -\R(\tau) S(\tau)I(\tau), ~~~~~~~~~~~~ \tau \in [0,T], &  \\
		&\dot I(\tau) = \R(\tau) S(\tau)I(\tau) - I(\tau), ~~~~~~ \tau \in [0,T], &  \\
		&(S(\tau),I(\tau)) \in \setX,~~~~~~~~~~~~~~~~~~~~~~~  \tau \in [0,T], &  \\
		&I(\tau) \leq I_{max},~~~~~~~~~~~~~~~~~~~~~~~~~~~~~~  \tau \in [0,T], & \\
		&S(T) = S^*,~~\R(\cdot) \in \Omega_\R, &
	\end{array}
\end{eqnarray}
where $T \! > \! \tau_f$ is a large enough (possibly infinite) time that covers the whole dynamic of the epidemic and $\Omega_\R$ is the set defined in equation \eqref{eq:Rform}. Conditions $I(\tau) \! \leq \! I_{max}$ and $S(T) \!=\! S^*$ force variable $I(\tau)$ to be smaller than an externally imposed maximum $I_{max}$ ($I_{max} \! \geq \! I(\tau_s)$), at every time $\tau \! \in \! [0,T]$, and variable $S(\tau)$ to be equal to $S^*$ at the end time $T$, respectively. In Problem $\mathcal P_{opt}(S(0),I(0),S^*,I_{max};\R(\cdot))$, $S(0)$, $I(0)$, $S^*$ and $I_{max}$ are optimization parameters, while function $\R(\cdot)$ is the optimization variable.

The next Theorem establishes that Problem $\mathcal P_{opt}$ is well-posed, achieves the epidemiological control objective, and gives better results, in general, than the two previous strategies in terms of the social/economic control objectives.
\begin{theo}\label{the:existsol}
	Consider system \eqref{eq:SIRnondim} with initial conditions $(S(0),I(0))=(1-\epsilon,\epsilon)$, $0 < \epsilon \ll 1$, and $\R(0)=\oR$ such that $S(0)\!>\!S^*$. Consider a given $I_{max}$ and consider also that $\R(\cdot) \in \Omega_{\R}$. Then, the solution of Problem $\mathcal P_{opt}(S(0),I(0),S^*,I_{max};\R(\cdot))$, denoted as $\R^{opt}$, produces $(S(\tau),I(\tau))$ for $\tau \! \geq \!0$ that arbitrarily approaches the epidemiological control objectives as $\tau_f \! \rightarrow \! \infty$ and, furthermore, minimizes the SDI, which constitutes the social/economic control objective.
\end{theo}
\begin{pf}
	See Appendix 2.
\end{pf}
\begin{exmp}
Figure \ref{fig:OptCont} shows the system evolution when the optimal social distancing $\R^{opt}$ is applied to the system under study. Time $T$ was selected to be $270$ days, and after this time, the system is left in open-loop. As it can be seen, a significant improvement is achieved in terms of the economic/social control objective: the SDI drop from $SDI\approx 300$ in the two previous strategies (goldilocks single interval and "wait, maintain, suspend" strategies) to $SDI \approx 192$, while the epidemiological control objectives are practically reached. 
As expected, the way the control problem is put allows us to account for both control objectives, without the need for competition between them, as is usually the case in other proposals. According to the time evolution of $\R(\cdot)$, the control seems to separate the epidemiological objectives over time (as it is done in the "wait, maintain, suspend" strategy): first, the control is devoted to controlling the IPP (from $\tau \approx$ $48$ to $\tau \approx 67$ days) and, once no further increments of $I$ can take place, it tries to reach the $S(\tau_f) \approx S^*$, in such a way that $S_\infty$ approaches $S^*$ at steady-state, and the EFS is minimized.
\begin{figure*}
	\centering
	\begin{subfigure}{.4\textwidth}
		\centering
		\includegraphics[width=1 \columnwidth]{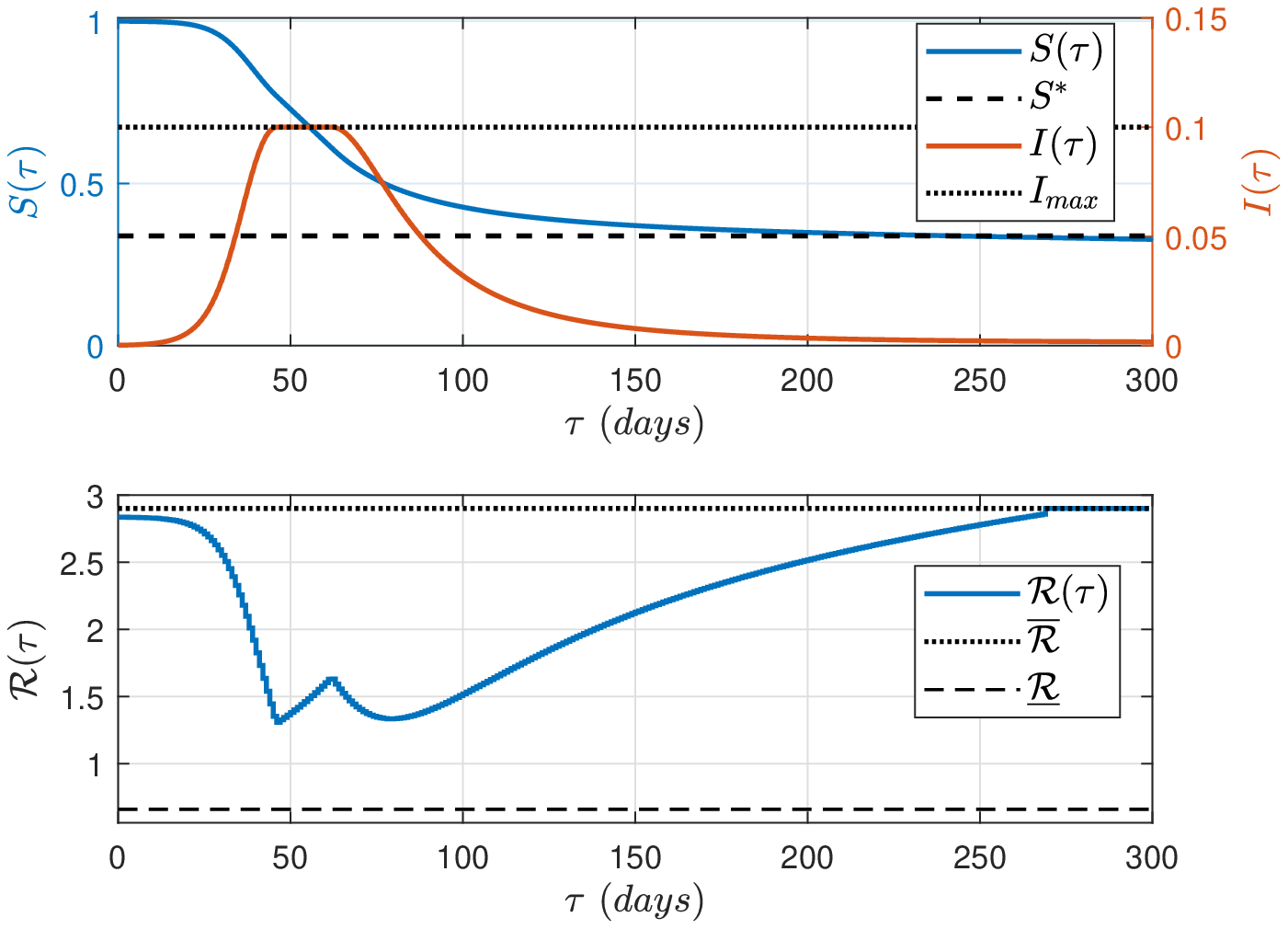}
	\end{subfigure}%
	\hspace*{0.2truecm}
	\begin{subfigure}{.4\textwidth}
		\centering
		\includegraphics[width=1 \columnwidth]{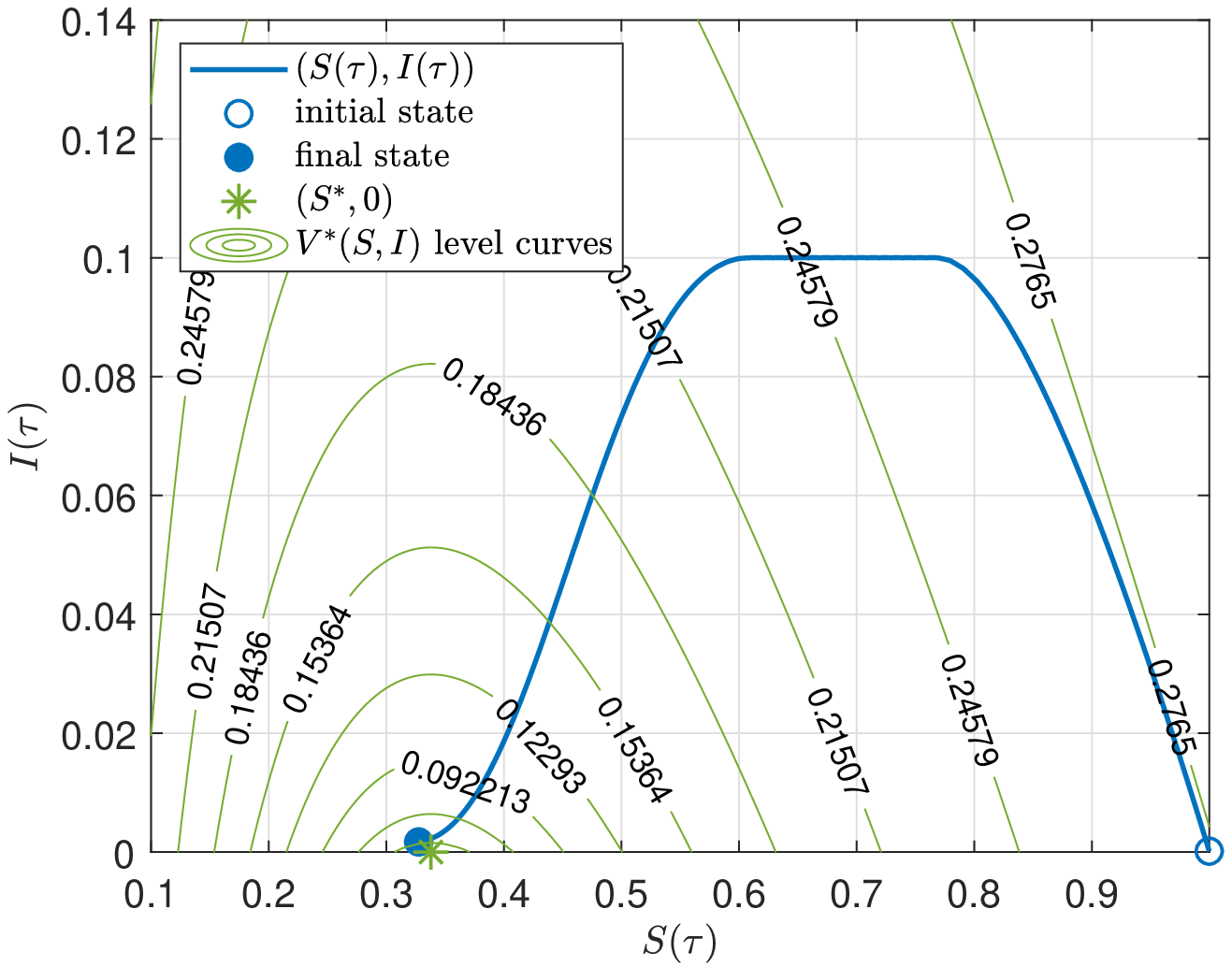}
	\end{subfigure}
	\caption{\small{$S,V$ and $\R$ time evolution, (left), and phase portrait, (right). System with $S^*=0.3448$ and $I_{max}=0.1$, when the control sequence $\R^{opt}$ is applied. The indexes are: EFS$=0.6725$, IPP$=0.1$, SDI$=192.4838$.}}
	\label{fig:OptCont}
\end{figure*}

Now, to emphasize the fact that problem $\mathcal P_{opt}$ is well-posed, we simulate the case in which the objective function $V$ is devoted to minimizing $\int_{0}^{T} I(t) dt$, as well as the social distancing severity. Let consider the objective function given by
\begin{eqnarray} \label{eq:valtern}
	V(\R(\cdot)) = \int_{0}^{T} \alpha_I I(t) + \alpha_\R(\oR -\R(t)) dt,
\end{eqnarray}
where $\alpha_I$ and $\alpha_\R$ are penalization weights that handle the relative importance of each term. Consider also that constraints $I(\tau) \leq I_{max}$, $\tau \in [0,T]$ and $S(T) = S^*$ are removed from Problem $\mathcal P_{opt}$. Figure \ref{fig:OptSCont} shows a simulation with $\alpha_I=1$ and $\alpha_\R$ taking the values: $0.25$ and $0.35$. As shown, there is not a combination of the penalties able to produce an optimal control strategy that achieves the epidemiological objective. For $\alpha_\R=0.25$, $\R(\cdot)$ is able to keep $I(\tau)$ under $I_{max}$, but $S(\tau_f)$ is two large, so a second epidemic wave is experienced, over the final simulation time. For $\alpha_\R=0.35$, on the other hand, $S_\infty \approx S^*$, but $I(\tau)$ overpass $I_{max}$, by a significant amount. 
\begin{figure*}
	\centering
	\begin{subfigure}{.4\textwidth}
		\centering
		\includegraphics[width=1 \columnwidth]{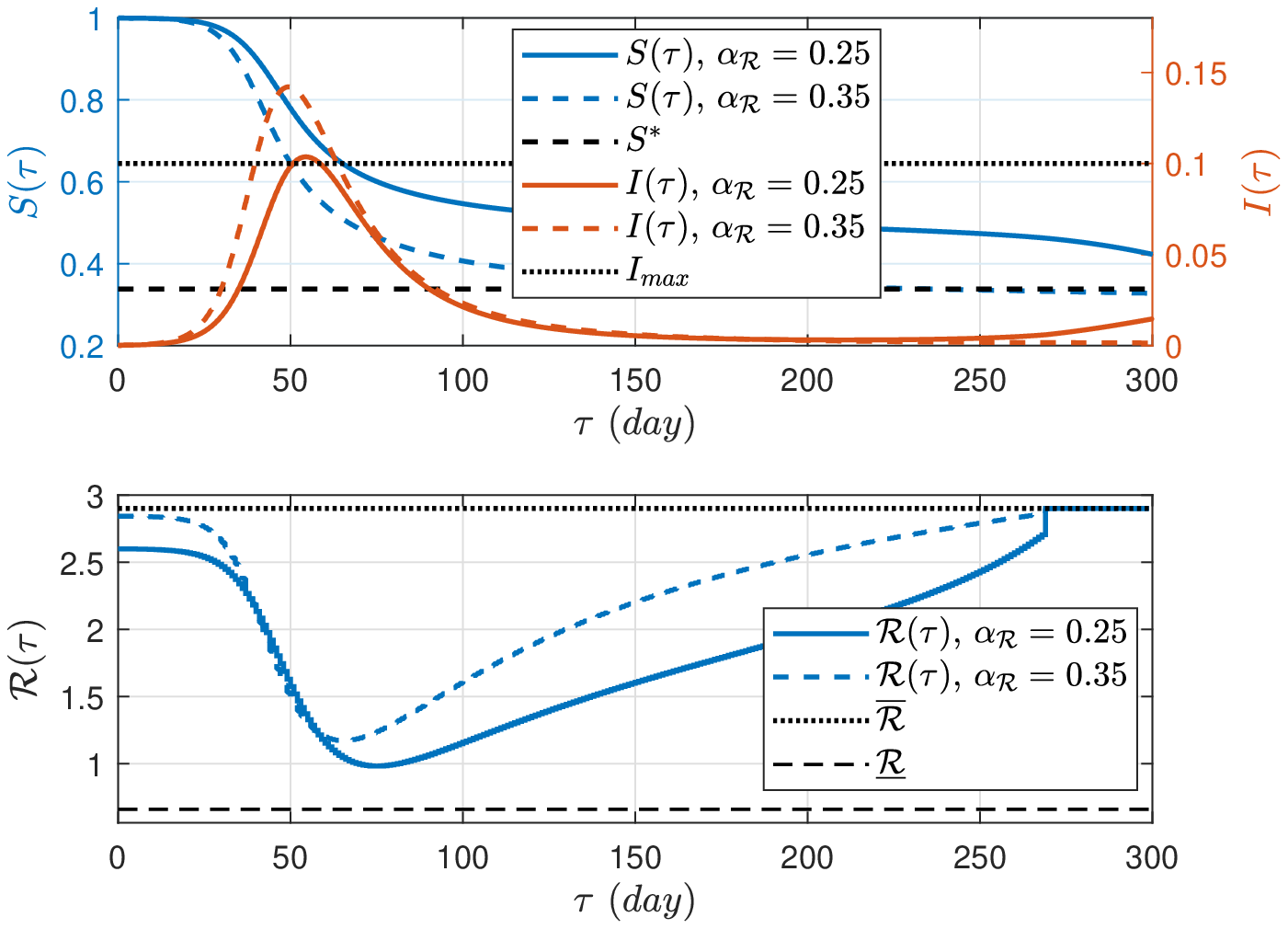}
	\end{subfigure}%
	\hspace*{0.2truecm}
	\begin{subfigure}{.4\textwidth}
		\centering
		\includegraphics[width=1 \columnwidth]{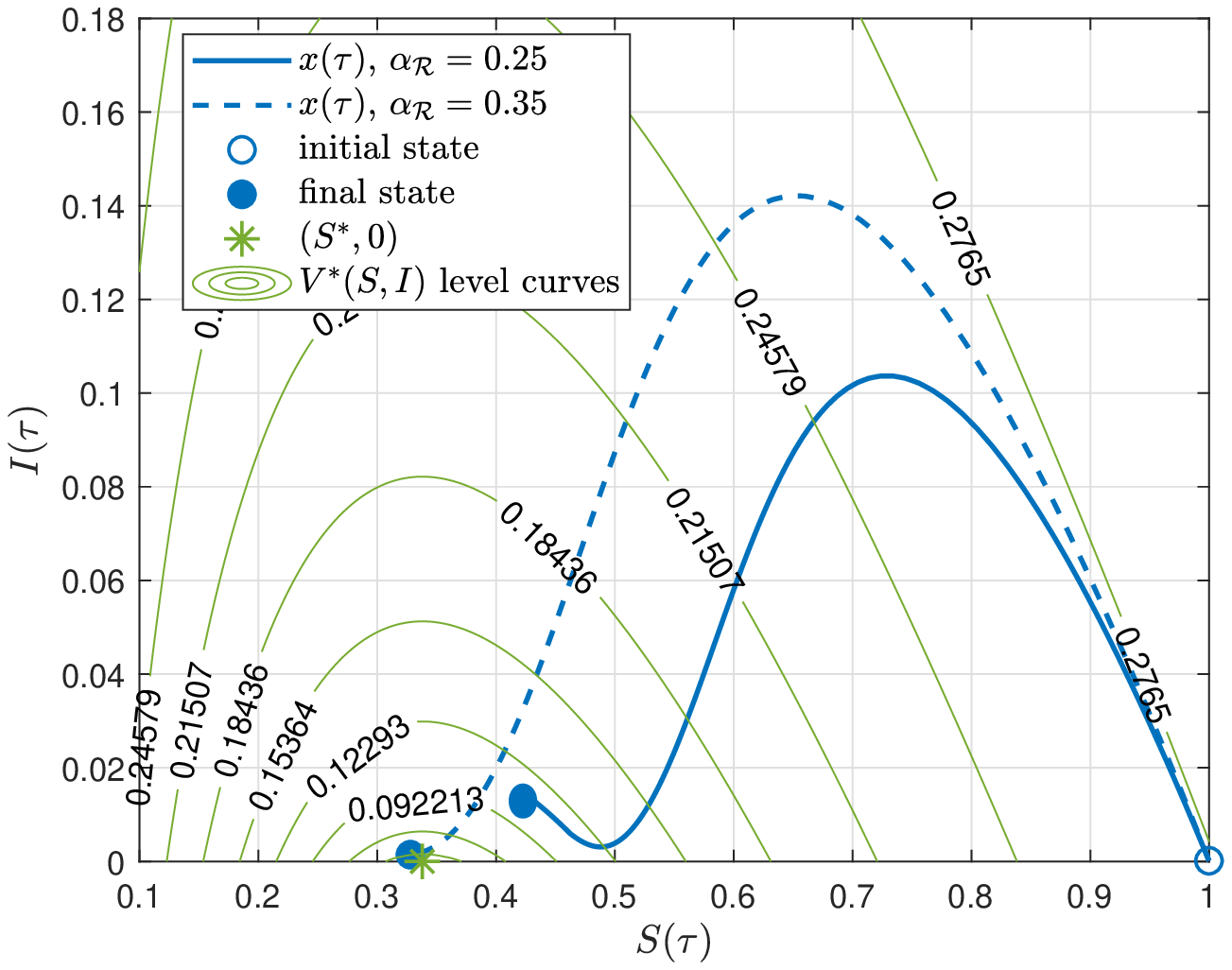}
	\end{subfigure}
	\caption{\small{$S,V$ and $\R$ time evolution, (left), and phase portrait, (right). System with $S^*=0.3448$ and $I_{max}=0.1$, when objective function \eqref{eq:valtern} is used for the optimal control problem $\mathcal P_{opt}$, and constraints $I(\tau) \leq I_{max}$, $\tau \in [0,T]$ and $S(T) = S^*$ are removed. There is no combination of $\alpha_I$ and $\alpha_\R$ that achieves the epidemiological control objectives.}}
	\label{fig:OptSCont}
\end{figure*}

From the point of view of realistic strategies, however, the optimal solution $\R^{opt}$ is still far from being directly applicable. First, $\R^{opt}(\cdot)$ varies continuously over time in some periods of time and, second, it takes any real value in $[\uR,\oR]$. Social distancing measures are usually described by only a few possible values, ranging from hard to soft (or even no) interventions and they are applied only for a certain fixed period of time, computed in weeks or months.
\end{exmp}
%
\subsection{Toward a more realistic  strategy}\label{sec:realcont} 
%
For the sake of completeness of the presentation, and taking advantage of the well-posed problem $P_{opt}$, let us assume that $\R(\cdot)$ can now take only four values in $[\uR,\oR]$, say, $[\uR,\R_1,\R_2,\oR] = [0.660, 1.407, 2.153, 2.900]$ and, furthermore, each of these value can be used for a period of time no smaller than ten days. These two extra conditions can be modelled by modifying the set $\Omega_\R$ in Problem $\mathcal P_{opt}$. 
\begin{exmp}
The simulation results corresponding to this new scenario are shown in Figure \ref{fig:MPCCont}, where it can be seen that the epidemiological objective is still achieved, but in a rather different way. The piece-wise control actions still seem to separate the epidemiological objectives, as before: the first part of the control action (from $\tau \approx 30$ to $\tau \approx 100$ days) is devoted to maintaining $I(\tau)$ under $I_{max}$, while the second one is devoted to steers $S$ to $S^*$, at steady-state. The SDI is given by $201.6$, which is higher than the optimal one, but still significantly better than the two first cases.
\begin{figure*}
	\centering
	\begin{subfigure}{.4\textwidth}
		\centering
		\includegraphics[width=1 \columnwidth]{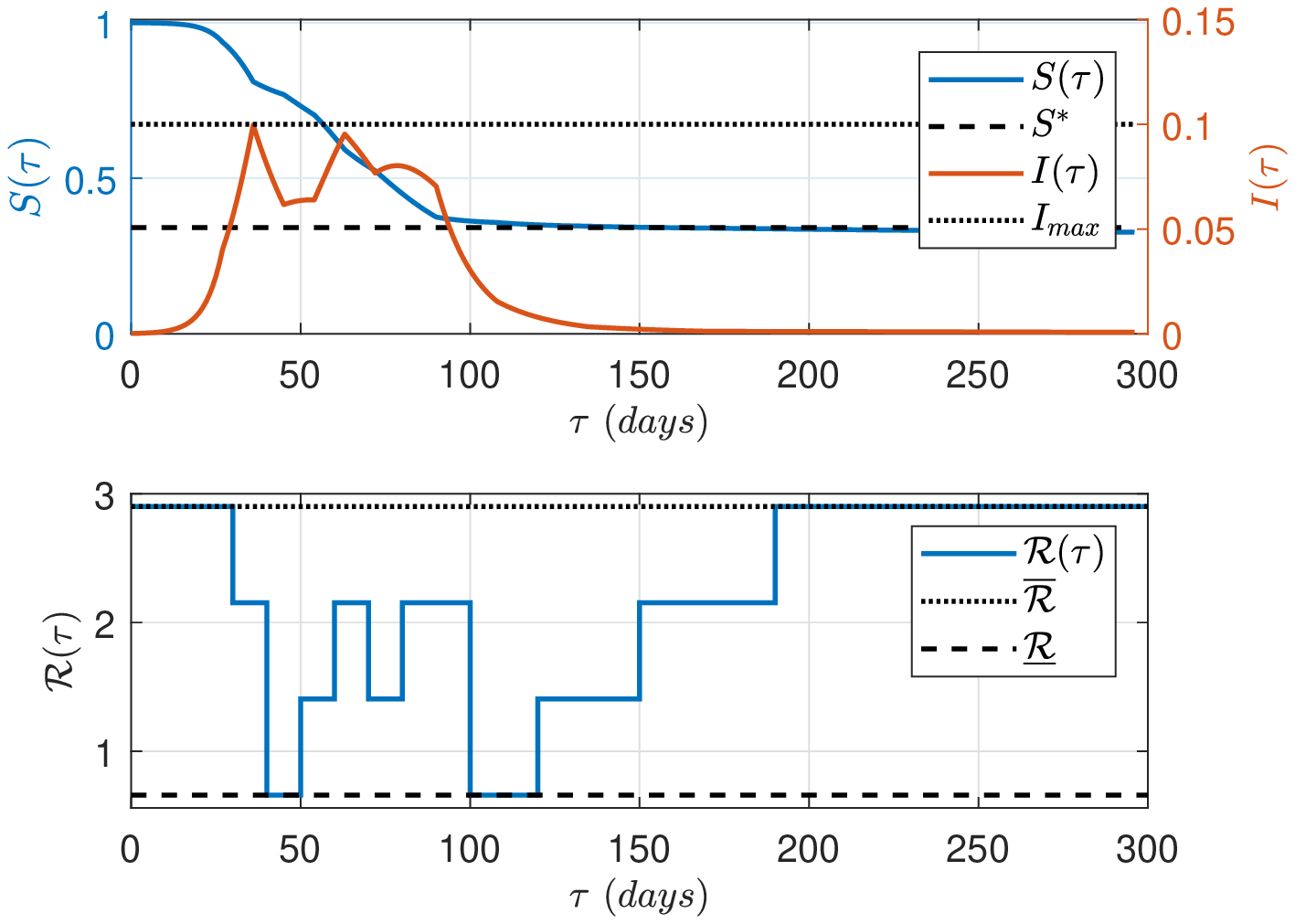}
	\end{subfigure}%
	\hspace*{0.2truecm}
	\begin{subfigure}{.4\textwidth}
		\centering
		\includegraphics[width=1 \columnwidth]{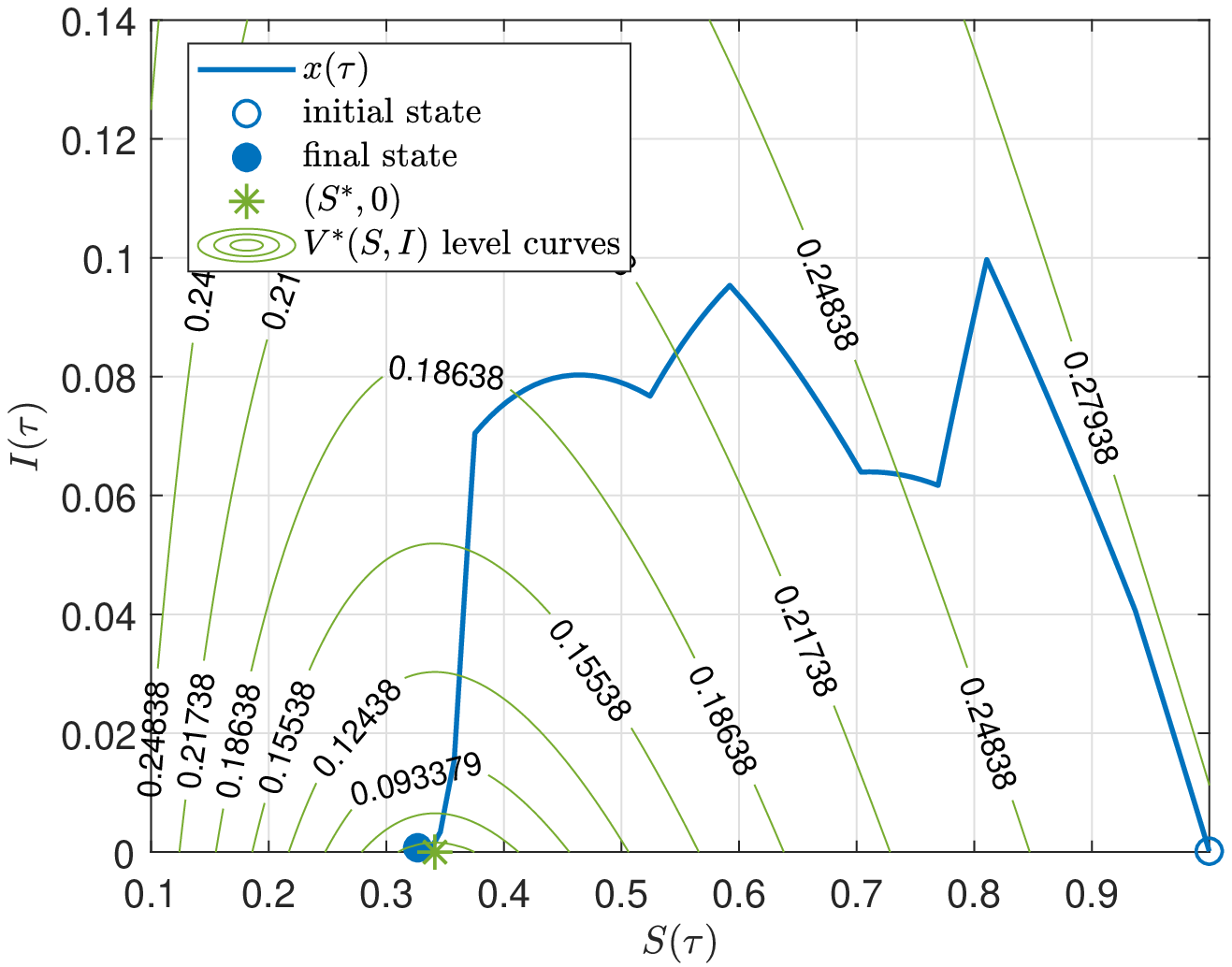}
	\end{subfigure}
	\caption{\small{$S,V$ and $\R$ time evolution, (left), and phase portrait, (right). System with $S^*=0.3448$ and $I_{max}=0.1$, when additional realistic conditions on $\R(\cdot)$ - as quantized levels and piece-wise constant behavior - are incorporated to the constraints of Problem $\mathcal P_{opt}$. Although the control is significantly different from the one obtained before, it still accounts for both control objectives. The indexes are: EFS$=0.6733$, IPP$=0.1$, SDI$=201.6$.}}
	\label{fig:MPCCont}
\end{figure*}
\end{exmp}
The main advantages of the proposed optimal control problem are: 
	(i) The cost function only penalizes the social distancing. So, once $S$ approaches $S^*$, for a time $\tau$ large enough, this minimization gradually suspends the social distancing, ensuring that social distancing will end at the smallest possible time. 
	(ii) The cost function does not consider $I(\tau)$, so there is no competition between the EFS and IPP objectives. 
	(iii) The constraint on the peak of $I$ ($I(\tau) \! \leq \! I_{max}$, $\tau \! \in \![\tau_s,\tau_f]$) allows us to arbitrarily control the IPP, by selecting arbitrarily small values of $I_{max}$. The effect of smaller values of $I_{max}$ will be prolonged social distancing, but it does not affect the optimality of the EFS.
	(iv) The degrees of freedom to select the form of $\R^{opt}(\cdot)$ allow one to take realistic social distancing measures: $\R^{opt}$ taking only a few possible values in $[\uR,\oR]$, corresponding to the degree of severity of the intervention, and lasting fixed period of time (of several days) before it is changed.

The main drawback, on the other hand, is that the optimization problem does not represent, in the current form, a proper feedback. It does not consider real-time information for computing the control action and, so, small differences between the model and the reality, or even an error in the estimation of $S$ and $I$, will necessarily produce performances far from the optimal one. 
%

\section{Conclusions and future works}

In this work, the equilibrium and stability of SIR-type models were characterized from a set-based perspective. Based on this characterization, the infected peak prevalence and the epidemic final size were studied, and concrete epidemiological objectives involving both indexes were proposed. It is shown that even with simple control strategies, a proper schedule for the social distancing can be designed to account for such an objective, at least from a theoretical point of view. If social and economic side effects of the social distancing are also considered - with leads to a non-trivial optimal control problem - it is shown that a solution can be found that, in addition to achieving the epidemiological objective, minimizes (as long as possible) the side effects.

Although theoretical, this study clearly shows what can and cannot be done to handle an epidemic: (i) just wait (with any control strategy) for reaching the herd immunity is not an option, since this way the EFS is not minimized (the herd immunity must be reached at steady-state); (ii) over-control the epidemic by means of hard interventions is not an option, since a fictitious steady-state will be reached that produces a second epidemic wave and, again, EFS will not be minimized; (iii) to directly minimize the IPP as a pure control objective gives not satisfactory results since scenarios far from the optimal one will be systematically achieved.

The next step toward a more realistic is to propose a proper feedback controller able to account for uncertainties/disturbances over short, updated, time horizons. Model Predictive Control (MPC) approach appears to be the right framework to account for such a challenge since they represent the best choice to close the loop under constraints. One step ahead is the explicit consideration of vaccination, immune system, and mutation strains, which also implies an extension of the model.

\section*{Appendix 1. Stability theory}\label{sec:app1}

All the following definitions are referred to system 
\begin{eqnarray}\label{eq:difeqini}
	\dot x (t) = f(x(t)),~~x(0)=x_0,
\end{eqnarray}
where $x$ is the system state constrained to be in $\X \subseteq \mathbb R^n$, $f$ is a Lipschitz continuous nonlinear function,
and $\phi(t;x)$ is the solution for time $t$ and initial condition $x$.
\begin{defn}[Equilibrium set]
	Consider system \ref{eq:difeqini} constrained by $\X$. The set $\setX_s \subset \X$ is an equilibrium set if each point $x \in \setX_s$ is such that $f(x) =  0$ (this implying that $\phi(t;x)  =  x$ for all $t \geq 0$).
\end{defn}
\begin{defn}[Attractivity of an equilibrium set]\label{def:attrac_set}
	Consider system \ref{eq:difeqini} constrained by $\X$ and a set $\setX \subseteq \X$. A closed equilibrium set $\setX_s \subset \setX$ is attractive in $\setX$	if $\lim_{t \rightarrow \infty} \|\phi(t;x)\|_{\setX_s} =0$ for all $x \in \setX$.
\end{defn}

A closed subset of an attractive set (for instance, a single equilibrium point) is not necessarily attractive. On the other hand, any set containing an attractive set is attractive, so the significant attractivity concept in a constrained system is given by the smallest one.

\begin{defn}[Local $\epsilon - \delta$ stability of an equilibrium set]\label{def:eps_del_stab}
	Consider system \ref{eq:difeqini} constrained by $\X$. A closed equilibrium set $\setX_s \subset \X$ is $\epsilon-\delta$ locally stable
	if for all $\epsilon >0$ it there exists $\delta>0$ such that in a given boundary of $\setX_s$, $\|x\|_{\setX_s} <\delta$, it
	follows that $ \|\phi(t;x)\|_{\setX_s} < \epsilon$, for all $t \geq 0$.
\end{defn}

Unlike attractive sets, a set containing a locally $\epsilon-\delta$ stable equilibrium set is not necessarily locally $\epsilon-\delta$ stable. Closed subsets of a locally $\epsilon-\delta$ stable equilibrium set is $\epsilon-\delta$ stable.

\begin{defn}[Asymptotic stability (AS) of an equilibrium set]\label{def:AS}
	Consider system \ref{eq:difeqini} constrained by $\X$ and a set $\setX \subseteq \X$. A closed equilibrium set $\setX_s \in \X$ is asymptotically stable (AS) in $\setX$ if it is $\epsilon-\delta$ locally stable and attractive in $\setX$.
\end{defn}

Next, the Lyapunov theorem, which refers to single equilibrium points and provides sufficient conditions for both, local $\epsilon-\delta$  stability and asymptotic stability, is introduced.
\begin{theo}[Lyapunov theorem (\cite{khalil2002nonlinear})]\label{theo:lyap}
	Consider system \ref{eq:difeqini} constrained by $\X$ and an equilibrium state $x_s \in \X$. Let consider a 
	function $V(x): \mathbb R^n \rightarrow \mathbb R$ such that $V(x)>0$ for $x \neq x_s$, $V(x_s)=0$ and $\dot{V}(x(t)) \leq 0$, denoted as Lyapunov function.
	Then, the existence of such a function in a boundary of $x_s$ implies that $x_s \in \setX_s$ is locally $\epsilon-\delta$ stable. If in addition
	$\dot{V}(x(t)) < 0$ for all $x\in \setX \subseteq \X$ such that $x \neq x_s$, and $\dot{V}(x_s) = 0$, then $x_s$ is asymptotically stable in $\setX$.
\end{theo}
%
\section*{Appendix 2. Technical proofs}\label{sec:app2}
%
\subsection*{Proof of Lemma \ref{lem:Sinf_opt}}
%
\begin{pf}
	Denote, for the sake of simplicity, $S(\tau_0)=S$, $I(\tau_0)=I$ and $\R=\oR$. Then, define $S_\infty^{op}(\delta):=\max_{S,I} \{S_\infty(\R,S,I) : (S,I) \in \setE(\delta) \}$,
	%
	%
	where $\setE(\delta):=\{(S,I) \in \mathbb R^2\!:\! S\!\in\! [0,1],~I \!\in\! [\delta,1]\}$ is a set of initial conditions with $I\!\geq\!\delta$, for some fixed $\delta \!\in\! [0,1]$, and $S^*$ is defined in \eqref{eq:Sstar}. Furthermore, define the maximizer initial conditions as $(S^{op}(\delta),I^{op}(\delta))\!\!:=\!\!arg\max_{S,I} \{S_\infty(\oR,S,I) \!:\! (S,I) \!\in\! \setE \}$.
	%
	%
	We will show that $S_\infty^{op}(\delta) = - W(- \R S^* e^{-\R (S^*+\delta)})/\R$, $(S^{op}(\delta),I^{op}(\delta)) = (S^*,\delta)$ and, particularly, that $S_\infty^{op}:= S_\infty^{op}(0)=S^*$.

	According to \eqref{eq:Sinf}, $S_\infty$ is given by $S_\infty(\R,S,I):= - W(-f(\R,S,I))/\R$,
	%
	%
	with $f(\R,S,I):=\R S e^{-\R (S+I)}$. Given that $-W(-x)$ is an increasing (injective) function of $x \in [0,1/e]$ and $\R$ is fixed, then $S_\infty(\R,S,I)$ achieves its maximum over $\setE(\delta)$ at the same values of $S$ and $I$ as $f(\R,S,I)$ (next it is shown that $f(\R,S,I) \in [0,1/e]$ for all $(S,I) \in \setE(\delta)$) and $\delta \in [0,1]$. Then, we focus the attention in finding the maximum (and the maximizing variables) of $f(\R,S,I)$. Let denotes the maximum of $f$ as $f^{op}(\delta):=\max_{S,I} \{f(S,I) : (S,I) \in \setE(\delta) \}$, while the maximizing variables are $S^{op}(\delta)$ and $I^{op}(\delta)$.

	Given that the maximum of $f$ occurs at the minimal values of $I$, let us consider, for simplicity, that $g(S,I) = I-\delta$, in such a way that we want to solve $(S^{op}(\delta),I^{op}(\delta))=arg\max_{S,I} \{f(\R,S,I) : g(S,I) \leq 0\}$ (we ignore the conditions $0 \leq S \leq 1$ and $I \leq 1$, but it is easy to see the no maximum is achieved at the boundaries of these constraints). Then $\bigtriangledown f =[\frac{\partial f}{\partial S},\frac{\partial f}{\partial I}]$ $ = [\R e^{-\R (S+I)}(1-\R S),\R^2 S e^{-\R (S+I)}]$ and $\bigtriangledown g =[\frac{\partial g}{\partial S},~\frac{\partial g}{\partial I}] \!=\! [0,1]$. Optimality conditions can be written as $\bigtriangledown f = \lambda \! \bigtriangledown \!g$, where $\lambda \in \mathbb R_{\geq0}$ is a Lagrange multiplier. Then, $[\R e^{-\R (S^{op}(\delta)+I^{op}(\delta))}(1-\R S^{op}(\delta))$, $ \R^2 S^{op}(\delta) e^{-\R (S^{op}(\delta)+I^{op}(\delta))}]=[0,~\lambda]$, which implies that
	$\R e^{-\R (S^{op}(\delta)+I^{op}(\delta))}(1-\R S^{op}(\delta))=0$ and $\R^2 S^{op}(\delta) e^{-\R (S^{op}(\delta)+I^{op}(\delta))} = \lambda$. Since $\R>0$, the first equality implies that $1-\R S^{op}(\delta) =0$, or $S^{op}(\delta) = \min\{1/\R\}=S^*$ (since $S^{op}(\delta)$ must be in $[0,1]$). This way, the second equality reads $\R^2 S^{*} e^{-\R (S^{*}+I^{op}(\delta))} = \lambda$, which is true for any value of $I^{op}(\delta)\in [\delta,1]$ and $\lambda>0$. As we know that larger values of $f$ are achieved for smaller values of $I$, then, $I^{op}(\delta)=\delta$.

	The maximum value of $S_\infty$ is then given by $S_\infty^{op}(\delta) = S_\infty(S^{op}(\delta),I^{op}(\delta))$, which reads $S_\infty^{op}(\delta) =$ $-W(-\R S^{op}(\delta) $ $ e^{-\R (S^{op}(\delta)+I^{op}(\delta))})/\R $ $= - W(-\R S^* e^{-\R (S^*+\delta)})/\R$.
	%
	%
	If $\R\!\!\geq\!\!1$, then $\R S^* \!\!=\!\!1$, and so $S_\infty^{op}(\delta)\!\!=\!\! - W(- e^{-\R (S^*+\delta)})/\R$. On the other hand, if $\R\!\!<\!\!1$, then $\R S^* \!\!=\!\!\R$ and $S_\infty^{op}(\delta)\!\!=\!\! - W(- \R e^{-\R (S^*+\delta)})/\R$. Particularly, $S_\infty^{op}:=S_\infty^{op}(0)= - W(- e^{-1})/\R = 1/\R = S^*$,
	if $\R\geq1$, and $S_\infty^{op}(0) =- W(-\R e^{-\R})/\R = \R/\R =1 =S^*$,
	if $\R<1$. This concludes the proof. $\square$
\end{pf}

\begin{figure}
	\centering
	\includegraphics[width=0.6\columnwidth]{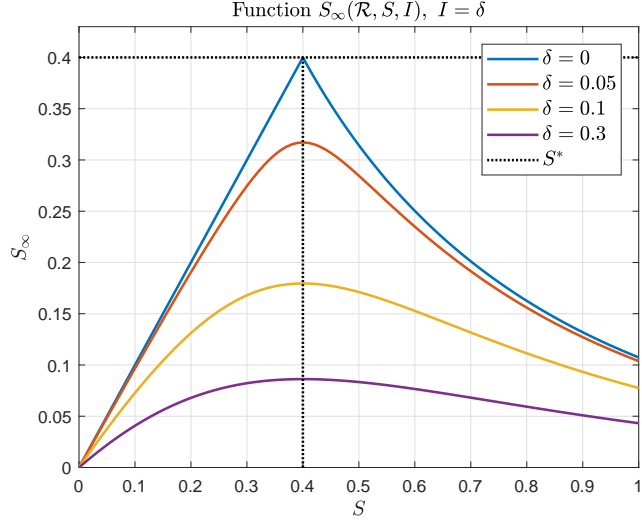}
	\caption{\small{Function $S_\infty(\R,S,I)$, with $\R=2.5$, for different values of $\delta$, when $S \in [0,1]$ and $I = \delta$. The maximum of $S_\infty(\R,S,\delta)$, for any fixed $\delta$, is achieved at $S=S^*$. Particularly, the maximum over $\delta \in [0,1]$ is achieved at $\delta=0$ and it is given by $S^*$.}}
	\label{fig:Sinf_S}
\end{figure}

In Lemma \ref{lem:Sinf_opt}, $\delta \in [0,1]$ represents the minimal initial value for $I$, in such a way that $\delta>0$ means that the initial conditions are not an equilibrium. Figure \ref{fig:Sinf_S} shows function $S_\infty(S,I)$ for different values of $\delta$, when $I =\delta$ and $\R=2.5$. Clearly, the maximum is achieved at $S=S^*$, when $\delta=0$. 
%
%
%
%
\subsection*{Proof of Theorem \ref{theo:stability}}
%
\begin{pf}
	The proof is divided into two parts. First it is shown that $\setX_s^{st}$ is the smallest attractive equilibrium set in $\setX$. Then, it is shown that $\setX_s^{st}$ is the largest locally $\epsilon-\delta$ stable equilibrium set in $\setX$ which, together with the previous results, implies that $\setX_s^{st}$ is the unique symptotically stable (AS) of system \eqref{eq:SIRnondim}, with a domain of attraction (DOA) given by $\setX$.

	\textit{Attractivity}:
	Consider equation \ref{eq:Sinf}.
	Since $W(z)$ is an increasing function (it goes from $-1$ at $z\!\!=\!\!-1/e$ to $0$ at $z\!\!=\!\!0$), it reaches its minimum at $z\!\!=\!\!-1/e$. $z(S(\tau_0),I(\tau_0))\!\!=\!\!-\R S(\tau_0)e^{-\R (S(\tau_0)\!+\!I(\tau_0))}$ reaches its maximum when $S(\tau_0)\!\!=\!\!S^*$, independently of the values of $\R$ and $I(\tau_0)$ (see Lemmas \ref{lem:Sinf_opt}), in which case it is $z(S(\tau_0),I(\tau_0))\!\!=\!\!1/e$. Then, $W(z)$ is bounded from above by $-1$, which means that $S^*$ is an upper bound for $S_\infty$. Therefore, $S_\infty \in [0,S^*]$, which shows the attractivity of $\setX_s^{st}$. 
	Figure \ref{fig:SinfFunc} shows a plot of $S_\infty$ as a function of $S(\tau_0)$ and $I(\tau_0)$ for a fixed value of $\R\!>\!1$ (a similar plot can be obtained for $\R\!<\!1$, in which $S_\infty$ reaches its maximum at the vertex of the domain, $S(\tau_0)\!\!=\!\!1$ and $I(\tau_0)\!\!=\!\!0$).

	To show that $\setX_s^{st}$ is the smallest attractive set in $\setX$, consider a state $(\bar S, \bar I):= \bar x \in \setX_s^{st}$ and an arbitrary small ball of radius $\epsilon > 0$, w.r.t. $\setX$, around it, $\mathbb{B}_{\epsilon}(\bar x) \in \setX$. Pick two arbitrary initial states $x_{0,1}=(S_{0,1},I_{0,1})$ and $x_{0,2}=(S_{0,2},I_{0,2})$ in $\mathbb{B}_{\epsilon}(\bar x)$, such that $S_{0,1} \neq S_{0,2}$. These two states converge, according to equation \eqref{eq:Sinf}, to $x_{\infty,1}=(S_{\infty,1},0)$ and $x_{\infty,2}=(S_{\infty,2},0)$, respectively, with $S_{\infty,1},S_{\infty,2} \in [0,S^*]$. Given that function $z(S(\tau_0),I(\tau_0))$ is monotone (injective) in $S(\tau_0)$ and $I(\tau_0)$, and $W(z)$ is monotone (injective) in $z$, then $S_{\infty,1} \neq S_{\infty,2}$. This means that, although both initial states converge to some state in $\setX_s^{st}$, they necessarily converge to different points. Therefore neither single states $\bar x \in \setX_s^{st}$ nor subsets of $\setX_s^{st}$ are attractive in $\setX$, which shows that $\setX_s^{st}$ is the smallest attractive set in $\setX$. 

	\textit{Local $\epsilon-\delta$ stability}:
	Let us consider a particular equilibrium point $\bar x := (\bar S,0)$, with $\bar S \in [0,S^*]$ (i.e., $\bar x \in \setX_s^{st}$).	Then a Lyapunov function candidate is given by:
	%
	\begin{eqnarray} \label{ec:lya1}
		V(x) := S- \bar S - \bar S \ln(\frac{S}{\bar S}) + I.
	\end{eqnarray}
	This function is continuous in $\setX$, is positive definite for all non-negative $x \neq \bar x$ and, furthermore, $V(\bar x)=0$. 
	Function $V$ evaluated at the solutions of system \eqref{eq:SIRnondim} reads:
	\begin{eqnarray} \label{ec:lya2}
		\frac{\partial V(x(\tau))}{\partial \tau} \!\! &=& \frac{\partial V}{\partial x} \dot{x}(\tau) 
		= \!\! \left[\frac{d V}{d S}~\frac{d V}{d I} \right] \left[
		\begin{array}{c}
			-\R S(\tau) I(\tau)  \\
			\R S(\tau) I(\tau)-I(t) 
		\end{array}\right]\nonumber\\
		&=& \!\! \left[(1-\frac{\bar S}{S(\tau)})~1 \right] \left[
		\begin{array}{c}
			-\R S(\tau) I(\tau)  \\
			\R S(\tau) I(\tau)-I(\tau) 
		\end{array}\right]\nonumber\\
		&=& \! I(\tau) (\R \bar S  - 1)
	\end{eqnarray}
	for $x(0)\in \setX$ and $\tau \geq 0$.
	Function $\dot{V}(x(\tau))$ depends on $x(\tau)$ only through $I(\tau)$. So, independently of the value of the parameter $\bar S$, $\dot{V}(x(\tau))=0$ for $I(\tau)\equiv 0$. This means that for any single $x(0) \in \setX_s$, $I(0)=0$ and so, $I(\tau)= 0$, for all $\tau\geq 0$. So $\dot{V}(x(\tau))$ is null for any $x(0) \in \setX_s$ (i.e, it is not only null for $x(0)=\bar x$ but for any $x(0) \in \setX_s$).

	On the other hand, for $x(0) \notin \setX_s$, function $\dot{V}(x(t))$ is negative, zero or positive, depending on if the parameter $\bar S$ is smaller, equal or greater than  $S^*= \min \{1,1/\R\}$, respectively, and this holds for all $x(0)\in \setX$ and $\tau \geq 0$. So, for any $\bar x \in \setX_s^{st}$, $\dot V(x(\tau))\leq 0$ (particularly, for $\bar x=(\bar S,0)=(S^*,0)$, $\dot V(x(\tau)) = 0$, for all $x(0) \in \setX$ and $\tau \geq 0$) which means that each $\bar x \in \setX_s^{st}$ is locally $\epsilon-\delta$ stable (see Theorem \ref{theo:lyap} in Appendix 1). Then, if every state in $\setX_s^{st}$ is locally $\epsilon-\delta$ stable, the whole set $\setX_s^{st}$ is locally $\epsilon-\delta$ stable.

	Finally, by following similar steps, it can be shown that $\setX_s^{un}$ is not $\epsilon-\delta$ stable, which implies that $\setX_s^{st}$ is also the largest locally $\epsilon-\delta$ stable set in $\setX$, which completes the proof. $\square$\\
\end{pf}
\begin{rem} \label{rem:nonas}
	In the previous proof, if we pick a particular $\bar x \in \setX_s^{st}$, then $\dot{V}(x(t))$ is not only null for $x(0)=\bar x$ but for all $x(0) \in \setX_s^{st}$, since in this case, $I(\tau)=0$, for $\tau\geq0$. This means that it is not true that $\dot{V}(x(t)) < 0$ for every $x \neq \bar x$, and this is the reason why single equilibrium points (and subsets of $\setX_s^{st}$) are $\epsilon-\delta$ stable, but not attractive. This is particularly true for the point $(S^*,0)$.
\end{rem}
%
%
\subsection*{Proof of Lemma \ref{lem:Suppbound}}
%
\begin{pf}
	The proof of (i) follows from Lemma \ref{lem:Sinf_opt}, by replacing $(S(\tau_0),I(\tau_0))$ by $(S(\tau_f),I(\tau_f))$, and the fact that the final intervention time, $\tau_f$, is finite. The proof of (ii) follows from Lemma \ref{lem:Sinf_opt}, the stability analysis made in Theorem \ref{theo:stability}, applied at $(S^*,0)$ (Corollary \ref{cor:Xst},(ii)), and Property \ref{propt:sinfty}, by replacing $(S(\tau_0),I(\tau_0))$ by $(S(\tau_f),I(\tau_f))$. $\square$
\end{pf}
%
%
\subsection*{Proof of Theorem \ref{theo:gold}}
%
\begin{pf}
	Consider that, for $S^*$ and $I_{max}$ coming from the Epidemiological Control Objectives, it there exists some $\tau_s$ for which $\R_{si}=\hat{\R}_{si}(I_{max},\tau_s) = \R_{si}^*(S^*,\tau_s)$. Then, by the definition of $\hat{\R}_{si}(I_{max},\tau_s)$, it follows that $\hat I(\hat{\R}_{si},S(\tau_s),I(\tau_s))=I_{max}$, which means that $I(\tau) \leq I_{max}$ for all $\tau>$. Furthermore, by the definition of $\R_{si}^*(S^*,\tau_s)$ it follows that, for a large enough $\tau_f$, $S(\tau_f)$ arbitrarily approaches $S^*$ while $I(\tau_f)$ approaches $0$. Then, by the stability results at $(S^*,0)$ (Corollary \ref{cor:Xst},(ii)), it follows that states $(S(\tau_f),I(\tau_f))$ arbitrarily close to $(S^*,0)$ (from above), produce states $(S(\tau),I(\tau))$ arbitrarily close to $(S^*,0)$ (from below), for $\tau >\tau_f$. Particularly, $S_\infty \approx S^*$. $\square$
\end{pf}
%
\subsection*{Proof of Theorem \ref{theo:mcsol}}
\begin{figure}
	\centering
	\includegraphics[width=0.6 \columnwidth]{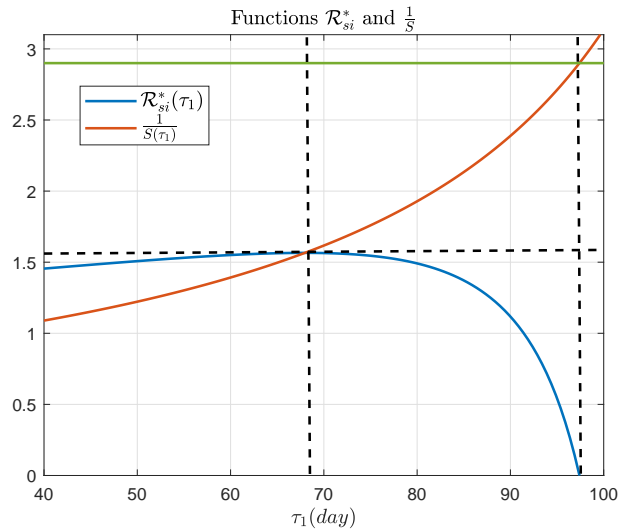}
	\caption{\small{Functions $\R_{si}^*(S^*,\tau_1)$ (blue), $\frac{1}{S(\tau_1)}$ (red) and $S^*$ (green), for $\oR=2.9$, $\tau_s=47.8$ and $I_{max}=0.1$. $\R_{si}^*(S^*,\tau_1)=\frac{1}{S(\tau_1)}$ at $\tau_1 = 68.7$ days, while $\R_{si}^*(S^*,\tau_1)=S^*$ at $\tau_1=97.4$ days.}}
	\label{fig:Rwms}
\end{figure}

%
\begin{pf}
	By hypothesis $I(\tau)$ starts growing from the small value $\epsilon$, at time $\tau=0$. Then, it there exists some $\tau_s>0$ such that $I(\tau_s)=I_{max}$. The control action $\R(\tau)=\frac{1}{S(\tau)}$ applied to system \eqref{eq:SIRnondim} for the period $[\tau_s,\tau_1)$ - that always fulfills the constraint $\R(\tau) \in [\uR,\oR]$, because $S(\tau_s)$ is close to, but small than, one - produces $\dot I(\tau) = 0$, which implies that $I(\tau)$ remains constant. Furthermore, $S(\tau)$ decreases linearly for $[\tau_s,\tau_1)$ (since $\dot S(\tau) = I(\tau) = I_{max}$). Now, time $\tau_1$ can be selected large enough for $I(\tau)$ not to increases for $\tau\geq \tau_1$, and small enough for $S(\tau)$ not decreases below $S^*$, for $\tau\geq \tau_1$. That is, the value of $S(\tau_1)$ must be smaller than the herd immunity value corresponding to system \eqref{eq:SIRnondim} with $\R(\tau)=\R^*_{si}$, but larger than the herd immunity value corresponding to $\oR$. This condition can be fulfilled if $\frac{1}{\R^{*}_{si}(S^{*},\tau_{1})} \geq S(\tau_{1}) \geq S^{*}$, or, the same
	$\R^*_{si}(S^*,\tau_1) \leq \frac{1}{S(\tau_1)} \leq \oR$,
	recalling that that $S(\tau) = (S(\tau_s)-I_{max} \tau_s) + I_{max} \tau$, for $\tau \in [\tau_s,\tau_1)$ and $\R^*_{si}(S^*,\tau_1)=\R^*_{si}(S^*,S(\tau_1),I(\tau_1))=\R^*_{si}(S^*,S(\tau_1),I_{max})$. Figure \ref{fig:Rwms} shows a plot of $\R_{si}^*(S^*,\tau_1)$ (blue), $\frac{1}{S(\tau_1)}$ (red) and $S^*$ (green).

	Any $\tau_{1}$ such that $\R^*_{si}(S^*,\tau_1) \leq \oR$, produces $S_\infty \approx S^*$ and $I(\tau)\leq I_{max}$. However, the best option, to minimize the severity of the social distancing is given by $\tau_1^*$; the time at which $\R^*_{si}(S^*,\tau_1^*) = \frac{1}{S(\tau_1^*)}$. $\square$
	%
	%
\end{pf}
%
\subsection*{Proof of Theorem \ref{the:existsol}}
%
\begin{pf}
	The proof follows form the fact that the "wait, maintain, suspend" strategy of Section \ref{sec:morcomplex} constitutes a feasible solution to problem $\mathcal P_{opt}(S(0),I(0),S^*,I_{max};\R(\cdot))$, which means that the optimal solution to this problem, $\R^{opt}(\cdot)$, will be one that produces, in general, a smaller cost $V(\R(\cdot))$ (and so, a smaller SDI) than the feasible one. Since, the "wait, maintain, suspend" strategy is not unique, problem $\mathcal P_{opt}$ will select a better one, and so it accounts for the economic/social control objective as well.
	$\square$\\
\end{pf}
\bibliographystyle{elsarticle-harv} 
\bibliography{Optimal_SIR_Biblio}
\end{document}